\documentclass{article}

\newcommand{\mK}{\mathbb{K}}

\usepackage{latexsym}
\usepackage{amsfonts}
\usepackage{amssymb}
\usepackage{cite}
\usepackage{color}

\newcommand{\brem}{\begin{remark}}
\newcommand{\erem}{\end{remark}}
\newcommand{\blem}{\begin{lemma}}
\newcommand{\elem}{\end{lemma}}
\newcommand{\bth}{\begin{theorem}}
\newcommand{\ethm}{\end{theorem}}
\newcommand{\benu}{\begin{enumerate}}
\newcommand{\eenu}{\end{enumerate}}
\newcommand{\bdes}{\begin{description}}
\newcommand{\edes}{\end{description}}
\newcommand{\bdf}{\begin{definition}}
\newcommand{\edf}{\end{definition}}
\newcommand{\bcor}{\begin{cor}}
\newcommand{\ecor}{\end{cor}}
\newcommand{\bprp}{\begin{proposition}}
\newcommand{\eprp}{\end{proposition}}
\newcommand{\bmlem}{\begin{mlemma}}
\newcommand{\emlem}{\end{mlemma}}
\newcommand{\bclm}{\begin{claim}}
\newcommand{\eclm}{\end{claim}}
\newcommand{\bprf}{{\bf Proof}.\hspace{2mm}}

\newcommand{\eprf}{\hspace*{\fill} $\Box$}

\newcommand{\beqn}{\begin{equation}}
\newcommand{\eeqn}{\end{equation}}
\newcommand{\beqnarr}{\begin{eqnarray}}
\newcommand{\eeqnarr}{\end{eqnarray}}
\newcommand{\beqnarrs}{\begin{eqnarray*}}
\newcommand{\eeqnarrs}{\end{eqnarray*}}

\newcommand{\spand}{\,\&\,}

\newtheorem{theorem}{Theorem}[section]
\newtheorem{definition}[theorem]{Definition}
\newtheorem{proposition}[theorem]{Proposition}
\newtheorem{lemma}[theorem]{Lemma}
\newtheorem{cor}[theorem]{Corollary}
\newtheorem{remark}[theorem]{Remark}
\newtheorem{mlemma}[theorem]{Main Lemma}
\newtheorem{claim}[theorem]{Claim}

\newcommand{\alp}{\alpha}

\newcommand{\veps}{\varepsilon}
\newcommand{\del}{\delta}
\newcommand{\Del}{\Delta}
\newcommand{\ome}{\omega}
\newcommand{\Ome}{\Omega}
\newcommand{\bet}{\beta}
\newcommand{\gam}{\gamma}
\newcommand{\Gam}{\Gamma}
\newcommand{\kap}{\kappa}
\newcommand{\sig}{\sigma}
\newcommand{\Sig}{\Sigma}
\newcommand{\tht}{\theta}
\newcommand{\Tht}{\Theta}
\newcommand{\lam}{\lambda}
\newcommand{\Lam}{\Lambda}
\newcommand{\vphi}{\varphi}

\newcommand{\fal}{\forall}
\newcommand{\exi}{\exists}

\newcommand{\Rarw }{\Rightarrow}

\newcommand{\lrarw}{\leftrightarrow}
\newcommand{\Lrarw}{\Leftrightarrow}

\newcommand{\cala}{{\cal A}}

\newcommand{\calc}{{\cal C}}

\newcommand{\calE}{{\cal E}}

\newcommand{\calg}{{\cal G}}
\newcommand{\calh}{{\cal H}}

\newcommand{\cals}{{\cal S}}

\newcommand{\calw}{{\cal W}}
\newcommand{\calx}{{\cal X}}
\newcommand{\caly}{{\cal Y}}

\newcommand{\la}{\langle}
\newcommand{\ra}{\rangle}

\newcommand{\rk}{\mbox{{\rm rk}}}

\newcommand{\sfk}{{\sf k}}

\newcommand{\setm}{\setminus}

\title{
A simplified ordinal analysis of first-order reflection
}
\author{Toshiyasu Arai
\\
Graduate School of Science,
Chiba University
\\
1-33, Yayoi-cho, Inage-ku,
Chiba, 263-8522, JAPAN
\\
tosarai@faculty.chiba-u.jp
\thanks{Current address:
Graduate School of Mathematical Sciences,
The University of Tokyo,
3-8-1 Komaba, Meguro-ku,
Tokyo 153-8914, JAPAN
tosarai@ms.u-tokyo.ac.jp
}
}
\date{}
\begin{document}
\maketitle

\begin{abstract}
In this note we give a simplified ordinal analysis of first-order reflection.
An ordinal notation system $OT$ is introduced based on $\psi$-functions.
Provable $\Sig_{1}$-sentences on $L_{\ome_{1}^{CK}}$ are bounded through cut-elimination
on operator controlled derivations. 
\end{abstract}

\section{Introduction}\label{sect:introduction}

Let $ORD$ denote the class of all ordinals, $A\subset ORD$ and $\alp$ a limit ordinal.
$\alp$ is said to be $\Pi_{n}$\textit{-reflecting} on $A$ iff
for any $\Pi_{n}$-formula $\phi(x)$ and any $b\in L_{\alp}$, if
$\la L_{\alp},\in\ra\models\phi(b)$,
then there exists a $\bet\in A\cap\alp$ such that $b\in L_{\bet}$ and
$\la L_{\bet},\in\ra\models\phi(b)$.
Let us write
$\alp\in rM_{n}(A):\Lrarw \alp \mbox{ is } \Pi_{n}\mbox{-reflecting on } A$.
Also $\alp$ {\rm is said to be} $\Pi_{n}$\textit{-reflecting} {\rm iff}
$\alp$ 
is $\Pi_{n}$-reflecting on $ORD$.

It is not hard for us to show that
the assumption that the universe is $\Pi_{n}$-reflecting
is proof-theoretically reducible to iterabilities  of
the lower operation $rM_{n-1}$ (and Mostowski collapsings), cf.\cite{consv}.

In this paper we aim an ordinal analysis of $\Pi_{n}$-reflection.
Such an analysis was done by Pohlers and Stegert \cite{Pohlers} using reflection configurations introduced in 
M. Rathjen\cite{RathjenAFML1},
and an alternative analysis in \cite{LMPS, WFnonmon2,PTPiN}
with the complicated combinatorial arguments of ordinal diagrams and finite proof figures.
Our approach is simpler in view of combinatorial arguments.
In \cite{LMPS}, a $\Pi_{n}$-reflecting universe is resolved in
ramified hierarchies of lower Mahlo operations, and ultimately in iterations of recursively Mahlo operations.
Our ramification process is akin to a tower, i.e., has an exponential structure.
It is natural that an exponential structure emerges in lowering and eliminating first-order formulas
(in reflections), cf.\,ordinal analysis for the fragments ${\rm I}\Sig_{n-3}$ of the first-order arithmetic.
Mahlo classes $Mh_{k}(\xi)$ defined in Definition \ref{df:Cpsiregularsm}
to resolve or approximate $\Pi_{n}$-reflection are based on similar structure.
As in Rathjen's analysis for $\Pi_{3}$-reflection in \cite{Rathjen94},
thinning operations are applied on the Mahlo classes $Mh_{k}(\xi)$, and this yields
an exponential structure similar to the one in \cite{LMPS} as follows.

Let us consider the simplest case $N=4$.
Let
$\Lam:=\veps_{\mK+1}$, the next epsilon number above the lease $\Pi_{4}$-reflecting ordinal $\mK$.
Roughly $\pi\in Mh_{3}(\xi)$ designates the fact that  an ordinal $\pi$
 is $\Pi_{3}$-reflecting on $Mh_{3}(\nu)$ for any $\nu<\xi<\Lam$. 
Suppose a $\Pi_{3}$-sentence $\tht$ on $L_{\pi}$ is derived from the assumption $\pi\in Mh_{3}(\xi)$.
We need to find an ordinal $\kap<\pi$ for which $L_{\kap}\models\tht$ holds.
It turns out that $\kap\in Mh_{2}(\Lam^{\xi}a)$ suffices for an ordinal $a<\Lam$, where
the ordinal $\kap$ in the class $Mh_{2}(\Lam^{\xi}a)$ is $\Pi_{2}$-reflecting
on classes $Mh_{2}(\Lam^{\xi}b)\cap Mh_{3}(\nu)$ for any $b<a$ and any $\nu<\xi$.
Note that the class $Mh_{2}(\Lam^{\xi}a)$ is not obtained through iterations of recursively Mahlo operations
since it involves $\Pi_{4}$-definable classes $Mh_{3}(\nu)$.
The classes $Mh_{3}(\nu)\,(\nu<\xi)$ for the assumption $\pi\in Mh_{3}(\xi)$ are thinned out
with the new classes $Mh_{2}(\Lam^{\xi}b)\,(b<\Lam)$, cf.\,Lemma \ref{lem:KppiNlower}.
\\

\noindent
Our theorem runs as follows.
Let $\mbox{{\sf KP}}\Pi_{N}$ denote the set theory for $\Pi_{N}$-reflecting universes, and
$\mbox{{\sf KP}}\ome$ the Kripke-Platek set theory with the axiom of infinity.
$OT$ is a computable notation system of ordinals defined in section \ref{subsec:decidable},
$\Ome=\ome_{1}^{CK}$ and $\psi_{\Ome}$ is a collapsing function such that $\psi_{\Ome}(\alp)<\Ome$.
$\mK$ is an ordinal term denoting the least $\Pi_{N}$-reflecting ordinal in the theorems. 

\begin{theorem}\label{thm:2}
Suppose ${\sf KP}\Pi_{N}\vdash\tht$
for a $\Sig_{1}(\Ome)$-sentence $\tht$.
Then we can find an $n<\ome$ such that for $\alp=\psi_{\Ome}(\ome_{n}(\mK+1))$,
$L_{\alp}\models\tht$.
\end{theorem}

Actually the bound is seen to be tight, cf.\,\cite{KPpiNwfprf}.
\bth\label{th:wf}
${\sf KP}\Pi_{N}$ proves that {\rm each} initial segment 
\\
$\{\alp\in OT: \alp<\psi_{\Ome}(\ome_{n}(\mK+1))\}\, (n=1,2,\ldots)$ is well-founded.
\end{theorem}

Thus the ordinal $\psi_{\Ome}(\veps_{\mK+1})$ is seen to be the proof-theoretic ordinal of ${\sf KP}\Pi_{N}$.

\begin{theorem}\label{th:main}
\[
\psi_{\Ome}(\veps_{\mK+1})=|{\sf KP}\Pi_{N}|_{\Sig_{1}^{\Ome}}
:=\min\{\alp\leq\ome_{1}^{CK}: \fal \tht\in\Sig_{1}({\sf KP}\Pi_{N}\vdash\tht^{L_{\Ome}} \Rarw L_{\alp}\models\tht)\}
.\]
\end{theorem}

$A\subset ORD$ is \textit{$\Pi^{1}_{n}$-indescribable} in $\alp$ iff for any $\Pi^{1}_{n}$-formula $\phi(X)$ and any $B\subset ORD$, 
if $\la L_{\alp},\in;B\cap\alp\ra\models \phi(B\cap\alp)$, then there exists a $\bet\in A\cap\alp$ such that
$\la L_{\bet},\in;B\cap\bet\ra\models\phi(B\cap\bet)$.
A regular cardinal $\pi$ is \textit{$\Pi^{1}_{n}$-indescribable} ifff $ORD$ is $\Pi^{1}_{n}$-indescribable in $\pi$.

Let us mention the contents of this paper. 
In the next section \ref{sect:ordinalnotation} we define simultaneously
iterated Skolem hulls $\mathcal{H}_{\alpha}(X)$ of sets $X$ of ordinals, ordinals 
$\psi^{\vec{\xi}}_{\kappa}(\alp)$ for regular cardinals 
$\kappa$, $\alp<\veps_{\mK+1}$ and sequences $\vec{\xi}=(\xi_{2},\ldots,\xi_{N-1})$
of ordinals $\xi_{i}<\veps_{\mK+2}$,
and classes $Mh^{\alp}_{k}(\xi)$ under the \textit{assumption} that a $\Pi^{1}_{N-2}$-indescribable cardinal $\mK$ exists.
It is shown that for $2\leq k<N$, $\alp<\veps_{\mK+1}$ and each $\xi<\veps_{\mK+2}$, 
$(\mK \mbox{ {\rm is a }} \Pi^{1}_{N-2}\mbox{{-indescribable cardinal}}) \to \mK\in Mh^{\alp}_{k}(\xi)$ in
${\sf ZF}+(V=L)$.

In section \ref{subsec:decidable} a computable notation system $OT$ of ordinals is extracted.
Following W. Buchholz\cite{Buchholz},
operator controlled derivations for $\mbox{{\sf KP}}\Pi_{N}$ is introduced in section \ref{sect:controlledOme},
 and inference rules for $\Pi_{N}$-reflection are eliminated from
derivations in section \ref{subsec:elimpi11}.
This completes a proof of Theorem \ref{thm:2} for an upper bound.
\\



IH denotes the Induction Hypothesis, MIH the Main IH and SIH the Subsidiary IH.
We are assuming tacitly the axiom of constructibility $V=L$.
Throughout of this paper $N\geq 3$ is a fixed integer.

\section{Ordinals for $\Pi_{N}$-reflection}\label{sect:ordinalnotation}
In this section 
we work in the set theory 
${\sf ZFLK}_{N}$ obtained from ${\sf ZFL}={\sf ZF}+(V=L)$
by adding the axiom $\exi \mK(\mK \mbox{ is } \Pi^{1}_{N-2}\mbox{-indescribable})$
for a fixed integer $N\geq 3$.
For ordinals $\alp$, $\veps(\alp)$ denotes the least epsilon number above $\alp$.

Let $ORD\subset V$ denote the class of ordinals, $\mK$ the least $\Pi^{1}_{N-2}$-indescribable cardinal, and
$Reg$ the set of regular ordinals below $\mK$.
$\Theta$ denotes finite sets of ordinals$\leq\mK$.
$u,v,w,x,y,z,\ldots$ range over sets in the universe,
$a,b,c,\alp,\bet,\gam,\ldots$ range over ordinals$<\Lam$,
$\xi,\zeta,\nu,\mu,\iota,\ldots$ range over ordinals$<\veps(\Lam)=\veps_{\mK+2}$,
$\vec{\xi},\vec{\zeta},\vec{\nu},\vec{\mu},\vec{\iota},\ldots$ range over finite sequences over ordinals$<\veps(\Lam)$,
and $\pi,\kap,\rho,\sig,\tau,\lam,\ldots$ range over regular ordinals.
$\tht$ denotes formulas.
\\

Let $\vec{\xi}=(\xi_{0},\ldots,\xi_{m-1})$ be a sequence of ordinals.
 The \textit{length} $lh(\vec{\xi}):=m$.
Sequences consisting of a single element $(\xi)$ is identified with the ordinal $\xi$,
and $\emptyset$ denotes the \textit{empty sequence}.
$\vec{0}$ denotes ambiguously a zero-sequence $(0,\ldots,0)$
with its length
$0\leq lh(\vec{0})\leq N-1$.
$\vec{\xi}*\vec{\mu}=(\xi_{0},\ldots,\xi_{m-1})*(\mu_{0},\ldots,\nu_{n-1})=(\xi_{0},\ldots,\xi_{m-1},\mu_{0},\ldots,\mu_{n-1})$
denotes the \textit{concatenated} sequence of $\vec{\xi}$ and $\vec{\mu}$.


$
\Lam=\veps(\mK)=\veps_{\mK+1}
$
denotes the next epsilon number above the least $\Pi_{N-2}$-
\\
indescribable cardinal
$\mK$, and
$\veps(\Lam)=\veps_{\mK+2}$ the next epsilon number above $\Lam$.

\bdf\label{df:Lam}
{\rm

For a non-zero ordinal $\xi<\veps(\Lam)$, its Cantor normal form with base $\Lam$
is uniquely determined as
\beqn\label{eq:CantornfLam}
\xi=_{NF}\sum_{i\leq m}\Lam^{\xi_{i}}a_{i}=\Lam^{\xi_{m}}a_{m}+\cdots+\Lam^{\xi_{0}}a_{0}
\eeqn
where
$\xi_{m}>\cdots>\xi_{0}, \, 0<a_{i}<\Lam$.

  \benu
  \item\label{df:Lam2.0}
  $K(\xi)=\{a_{i}:i\leq m\}\cup\bigcup\{K(\xi_{i}):i\leq m\}$ is the set of \textit{components} of $\xi$
  with $K(0)=\emptyset$.
For a sequence $\vec{\xi}=(\xi_{0},\ldots,\xi_{n-1})$ of ordinals $\xi_{i}<\veps(\Lam)$,
$K(\vec{\xi})
:=\bigcup\{K(\xi_{i}): i<n\}$.

 \item\label{df:Exp2.1}
 For $\xi>1$,
 $te(\xi)=\xi_{0}$ in (\ref{eq:CantornfLam}) is the \textit{tail exponent}, 
 and
 $he(\xi)=\xi_{m}$ is the \textit{head exponent} of $\xi$, resp.
The \textit{head} $Hd(\xi):=\Lam^{\xi_{m}}a_{m}$, and the \textit{tail}
 $Tl(\xi):=\Lam^{\xi_{0}}a_{0}$ of $\xi$.


 \item\label{df:Exp2.2h}
 $he^{(i)}(\xi)$ is the \textit{$i$-th head exponent} of $\xi$, defined recursively by
 \\
 $he^{(0)}(\xi)=\xi$, $he^{(i+1)}(\xi)=he(he^{(i)}(\xi))$.
 
 The \textit{$i$-th tail exponent} $te^{(i)}(\xi)$ is defined similarly.
 
 \item\label{df:Exp2.3}
 $\zeta$ is a \textit{part} of $\xi$, denoted by
 $\zeta\leq_{pt}\xi$ iff
 \\
 $\zeta=_{NF}\sum_{i\geq n}\Lam^{\xi_{i}}a_{i}=\Lam^{\xi_{m}}a_{m}+\cdots+\Lam^{\xi_{n}}a_{n}$ {\rm for an} $n\, (0\leq n\leq m+1)$.
 
$\zeta<_{pt}\xi:\Lrarw\zeta\leq_{pt}\xi \spand \zeta\neq\xi$.
 
\item\label{df:Exp2.4}
A sequence $\vec{\mu}=(\mu_{0},\ldots,\mu_{n})$ is an \textit{iterated tail parts} of $\xi$, denoted by
$\vec{\mu}\subset_{pt}\xi$
iff $ \mu_{0}\leq_{pt}\xi \spand \fal i<n(\mu_{i+1}\leq_{pt} te(\mu_{i}))$.

\item\label{df:Exp2.5}
$\vec{\nu}=(\nu_{0},\ldots,\nu_{n})*\vec{0}<\xi$ iff there exists a sequence
$\vec{\mu}=(\mu_{0},\ldots,\mu_{n})$ such that
$\vec{\mu}\subset_{pt}\xi$ and $\nu_{i}<\mu_{i}$ for every $i\leq n$.

\item\label{df:Exp2.8}
Let $\vec{\nu}=(\nu_{0},\ldots,\nu_{n})$ and $\vec{\xi}=(\xi_{0},\ldots,\xi_{n})$ be
sequences of ordinals in the same length, and $0\leq k\leq n$.

$\vec{\nu}<_{k}\vec{\xi}:\Lrarw \fal i<k(\nu_{i}\leq\xi_{i})\land 
(\nu_{k},\ldots,\nu_{n})< \xi_{k}$.

\item\label{df:Exp2.9}
 $\zeta$ is a \textit{step-down} of $\xi$, denoted by
 $\zeta<_{sd}\xi$ iff
 \\
 $\zeta=\Lam^{\xi_{m}}a_{m}+\cdots+\Lam^{\xi_{1}}a_{1}+\Lam^{\xi_{0}}b+\nu$ for some ordinals
 $b<a_{0}$ and $\nu<\Lam^{\xi_{0}}$.

\item\label{df:Exp2.10}
$\vec{\nu}=(\nu_{0},\ldots,\nu_{n})*\vec{0}<_{sd}\xi$ iff 
$\nu_{i}<_{sd}te^{(i)}(\xi)$ for every $i\leq n$.

\item\label{df:Exp2.11}
$\zeta\leq_{sp}\xi:\Lrarw\exi\mu\leq_{pt}\xi(\zeta\leq_{sd}\mu)$, and
$\zeta<_{sp}\xi:\Lrarw\exi\mu\leq_{pt}\xi(\zeta<_{sd}\mu)$.

\item\label{df:Exp2.12}
$\vec{\nu}<_{sp}\xi$ iff 
$\vec{\nu}<_{sd}\mu$ for a $\mu\leq_{pt}\xi$.

Let $p(\vec{\nu},\xi)$ denote the number $p\, (0\leq p<m)$ such that
$\xi=_{NF}\mu+\sum_{i<p}\Lambda^{\xi_{i}}a_{i}$
for
$\mu=\Lambda^{\xi_{m}}a_{m}+\cdots+\Lambda^{\xi_{p}}a_{p}$ and $\vec{\nu}<_{sd}\mu$.

\eenu

}
\edf

Note that $(\nu)*\vec{0}<\xi\Lrarw \nu<\xi$, and
$(\xi, te(\xi),te^{(2)}(\xi),\ldots)\subset_{pt}\xi$.
Also $\zeta<_{sd}\xi\Lrarw\zeta<\xi$ if $\xi<\Lam$.

\bprp\label{prp:headcomparison}
$\xi<\mu<\veps(\Lam)  \Rarw te(\xi)\leq he(\xi)\leq he(\mu)$.
\eprp


\bprp\label{prp:idless}
$\vec{\nu}<\xi\leq\zeta \Rarw \vec{\nu}<\zeta$.
\eprp
{\bf Proof}
by induction on the lengths $n=lh(\vec{\nu})$.
Let $\vec{\mu}=(\mu_{0},\ldots,\mu_{n-1})$ be a sequence for 
$\vec{\nu}=(\nu_{0},\ldots,\nu_{n-1})$ such that
$\vec{\mu}\subset_{pt}\xi$ and $\fal i\leq n-1(\nu_{i}<\mu_{i})$, 
cf.\,Definition \ref{df:Lam}.\ref{df:Exp2.5}.

If $n=1$, then $\nu_{0}<\mu_{0}\leq_{pt}\xi\leq\zeta$.
$\nu_{0}<\zeta\leq_{pt}\zeta$ yields $\vec{\nu}=(\nu_{0})<\zeta$.

Let $n>1$. We have $(\nu_{1},\ldots,\nu_{n-1})<te(\mu_{0})$ with $(\mu_{1},\ldots,\mu_{n-1})\subset_{pt}te(\mu_{0})$.
We show the existence of a $\lam$ such that $\mu_{0}\leq\lam\leq_{pt}\zeta$ and $te(\mu_{0})\leq te(\lam)$.
Then IH yields $(\nu_{1},\ldots,\nu_{n-1})<te(\lam)$, and $\vec{\nu}<\zeta$ follows.

If $\mu_{0}\leq_{pt}\zeta$, then $\lam=\mu_{0}$ works.
Suppose $\mu_{0}\not\leq_{pt}\zeta$.
On the other hand we have $\mu_{0}\leq_{pt}\xi\leq\zeta$.
This means that $\xi<\zeta$ and there exists a $\lam\leq_{pt}\zeta$ such that $\mu_{0}<\lam$ and $te(\mu_{0})\leq te(\lam)$.
\eprf



\subsection{Ordinals}


\bdf
{\rm
\benu
\item
For $i<\ome$ and $\xi<\veps(\Lam)$, 
$\Lam_{i}(\xi)$ is defined recursively by $\Lam_{0}(\xi)=\xi$ and
 $\Lam_{i+1}(\xi)=\Lam^{\Lam_{i}(\xi)}$.

\item\label{df:Lam3}
For $A\subset ORD$, limit ordinals $\alp$ and $i\geq 0$, let
$
\alp\in M_{2+i}(A)$ iff $A\cap\alp \mbox{ {\rm is} } \Pi^{1}_{i}\mbox{{\rm -indescribable in} } \alp
$.

 \item\label{df:Lam4}
 $\kap^{+}$ denotes the next regular ordinal above $\kap$.
 
 \item
 $\Ome_{\alp}:=\ome_{\alp}$ for $\alp>0$, $\Ome_{0}:=0$, and
 $\Ome=\Ome_{1}$.
\eenu
}
\edf

Define simultaneously 
classes $\mathcal{H}_{\alpha}(X)$, $Mh^{\alp}_{k}(\xi)$, and 
ordinals $\psi_{\kappa}^{\vec{\xi}}(\alp)$ as follows.
We see that these are 
$\Sig_{1}$-definable as a fixed point in ${\sf ZFL}$,
 cf. Proposition \ref{prp:definability}.

Let $a<\Lam$, and $\vphi$ denote  the binary Veblen function.
Let us define a Skolem hull $\mathcal{H}_{a}(X)$ of $\{0,\mK\}\cup X$ under the functions
$+,
 \alpha\mapsto\omega^{\alpha}, 
 (\alp,\bet)\mapsto\vphi\alp\bet\, (\alp,\bet<\mK),
 \alp\mapsto\Ome_{\alp}\, (\alp<\mK)$
 and $\psi$-functions.
$Reg$ denotes the set of regular ordinals$\leq\mK$.

\bdf\label{df:Cpsiregularsm}
{\rm
$\calh_{a}[Y](X):=\calh_{a}(Y\cup X)
$ for sets $Y\subset \mK$.

\benu
\item\label{df:Cpsiregularsm.1}
{\rm (Inductive definition of} $\mathcal{H}_{a}(X)${\rm ).}
\benu
\item\label{df:Cpsiregularsm.10}
$\{0,\mK\}\cup X\subset\mathcal{H}_{a}(X)$.

\item\label{df:Cpsiregularsm.11}
 $x, y \in \mathcal{H}_{a}(X) \Rightarrow x+ y\in \mathcal{H}_{a}(X)$,
 $x\in\calh_{a}(X)\Rarw \omega^{x}\in \mathcal{H}_{a}(X)$, {\rm and}
  $x, y \in \mathcal{H}_{a}(X)\cap\mK \Rightarrow \vphi x y\in \mathcal{H}_{a}(X)$.
 
\item\label{df:Cpsiregularsm.12}
$\mK>\alp\in\mathcal{H}_{a}(X) \Rightarrow \Ome_{\alp}\in\mathcal{H}_{a}(X)$.



\item\label{df:Cpsiregularsm.13}
If $\pi\in\mathcal{H}_{a}(X)\cap Reg$ and
$b\in\mathcal{H}_{a}(X)\cap a$, then
$\psi_{\pi}(b)\in\mathcal{H}_{a}(X)$.

\item\label{df:Cpsiregularsm.14}
If $\{b,\xi\}\subset\mathcal{H}_{a}(X)$ with $\xi\leq b<a$, then
$\kap=\psi_{\mK}^{\vec{0}*(\xi)}(b)\in\mathcal{H}_{a}(X)$, where
$lh(\vec{0})=N-3$.

\item\label{df:Cpsiregularsm.15}
Let 
$\{\pi,b,c\}\subset \calh_{a}(X)$ with $\pi<\mK$, $2\leq k<N-1$ an integer,
and $\vec{\xi}=(\xi_{2},\ldots,\xi_{k},\xi_{k+1})*\vec{0}$ a sequence of ordinals $\xi_{i}<\veps(\Lam)$ with $lh(\vec{0})=N-2-k$
such that $\xi_{k+1}\neq 0$ and
$K(\vec{\xi})\subset\calh_{a}(X)$.
Assume
$\max(K(\vec{\xi})\cup\{c\}) \leq b<a$,
and
$\pi\in Mh_{2}^{b}(\vec{\xi})$.

Then
$\kap=\psi_{\pi}^{\vec{\nu}}(b)\in\calh_{a}(X)$ for the sequence
$\vec{\nu}=(\xi_{2},\ldots,\xi_{k}+\Lam^{\xi_{k+1}}c)*\vec{0}$ with $lh(\vec{0})=N-1-k$.

\item\label{df:Cpsiregularsm.16}
Let 
$\{\pi,b\}\subset \calh_{a}(X)$ with $\pi<\mK$, and $0\neq\xi<\veps(\Lam)$ an ordinal
with $K(\xi)\subset\calh_{a}(X)$.
Let $\vec{\nu}=(\nu_{2},\ldots,\nu_{N-1})$ be a sequence of ordinals$<\veps(\Lam)$
such that $K(\vec{\nu})\subset\calh_{a}(X)$.
Assume $\max K(\vec{\nu}) \leq b<a$, $K(\vec{\nu})\subset\calh_{b}(\pi)$,
$\pi\in Mh_{2}^{b}(\xi)$, and
$\vec{\nu}<\xi$, cf.\,Definition \ref{df:Lam}.\ref{df:Exp2.5}.

Then
$\kap=\psi_{\pi}^{\vec{\nu}}(b)\in\calh_{a}(X)$.

\eenu

\item\label{df:Cpsiregularsm.2}
 (Definitions of $Mh^{a}_{k}(\xi)$ and $Mh^{a}_{k}(\vec{\xi})$)
\\
First let
$\mK\in Mh^{a}_{N}(0) :\Lrarw 
\mK\in M_{N} 
\Lrarw \mK \mbox{ {\rm is }} \Pi^{1}_{N-2}\mbox{{\rm -indescribable}}
$.

The classes $Mh^{a}_{k}(\xi)$ are defined for $2\leq k< N$,
and ordinals $a<\Lam$, $\xi<\veps(\Lam)$.
Let $\pi$ be a regular ordinal$\leq\mK$. Then for $\xi>0$
\beqnarr
&&
 \pi\in Mh^{a}_{k}(\xi)
  :\Lrarw 
 \{\pi,a\}\cup K(\xi)\subset\calh_{a}(\pi) \spand
\label{eq:dfMhkh}
\\
&&
 \fal\vec{\nu}<\xi \left(
 K(\vec{\nu})\subset\calh_{a}(\pi) 
 \Rarw \pi\in M_{k}(Mh^{a}_{k}(\vec{\nu}))\right)
\nonumber
\eeqnarr

 where $\vec{\nu}=(\nu_{k},\ldots,\nu_{n})\,(2\leq k\leq n\leq N-1)$ varies through 
 non-empty sequences of ordinals$<\veps(\Lam)$
and
\[
\pi\in Mh^{a}_{k}(\vec{\nu}):\Lrarw \pi\in\bigcap_{k\leq i\leq n}Mh^{a}_{i}(\nu_{i})
.\]
By convention, let for $2\leq k<N$,
$\pi\in Mh^{a}_{k}(0):\Lrarw \pi\in Mh^{a}_{2}(\emptyset):\Lrarw 
\pi \mbox{ {\rm is a limit ordinal}}
$.
Note that by letting $\vec{\nu}=(0)$,
$\pi\in Mh^{a}_{k}(\xi) \Rarw \pi\in M_{k}$ for $\xi>0$.
Also $\vec{0}<1$, and $Mh^{a}_{k}(1)=M_{k}$.

\item\label{df:Cpsiregularsm.3}
 (Definition of $\psi_{\pi}^{\vec{\xi}}(a)$)
\\
 Let $a<\Lam$ be an ordinal, $\pi\leq\mK$ a regular ordinal and
 $\vec{\xi}$ a sequence of ordinals$<\veps(\Lam)$ such that
 $lh(\vec{\xi})=N-2$.
Then let
{\small
\beqn\label{eq:Psivec}
\psi_{\pi}^{\vec{\xi}}(a)
 :=  \min(\{\pi\}\cup\{\kap\in Mh^{a}_{2}(\vec{\xi})\cap\pi:   \calh_{a}(\kap)\cap\pi\subset\kap ,
   K(\vec{\xi})\cup\{\pi,a\}\subset\calh_{a}(\kap)
\})
\eeqn
}
Let
$\psi_{\pi}a:=\psi_{\pi}^{\vec{0}}a$,
where $lh(\vec{0})=N-2$, $Mh^{a}_{2}(\vec{0})=Lim$, and $\pi\in M_{2}$, i.e.,
$\pi$ is a regular ordinal.

\eenu
}

\edf

Note that 
$\pi\in Mh^{a}_{k}(\xi) \Rarw \fal\nu<\xi \left(\pi\in M_{k}(Mh^{a}_{k}(\nu))\right)$,
since $(\nu)<\xi$ holds with $(\xi)\subset_{pt}\xi$ for $\nu<\xi$.

\bprp\label{prp:Cantorclosed}
$b+c\in\calh_{a}[\Tht](d) \Rarw c\in\calh_{a}[\Tht](d)$, and 
$\ome^{c}\in\calh_{a}[\Tht](d) \Rarw c\in\calh_{a}[\Tht](d)$.
\eprp

The following Proposition \ref{prp:definability} is easy to see.

\begin{proposition}\label{prp:definability}
Each of 
$x=\mathcal{H}_{a}(y)\, (a<\Lam,y<\mK)$,
$x=\psi_{\kappa}a$,
$x\in Mh^{a}_{k}(\xi)$ and $x=\psi^{\vec{\xi}}_{\kap}(a)$, 
is a $\Sig_{1}$-predicate as fixed points in $\mbox{{\sf ZFL}}$.
\end{proposition}
\bprf 
This is seen from the facts that there exists a universal $\Pi^{1}_{n}$-formula, and by using it,
$\alp\in M_{n}(x)$ iff $\la L_{\alp},\in\ra\models m_{n}(x\cap L_{\alp})$
for some $\Pi^{1}_{n+1}$-formula $m_{n}(R)$ with a unary predicate $R$.
\eprf
\\





Let $A(a)$ denote the conjunction of 
$\forall u<\mK \exists ! x[x=\mathcal{H}_{a}(u)]$, and
\\
$\fal \vec{\xi}\fal x(
\max K(\vec{\xi})\leq a 
\spand K(\vec{\xi})\cup\{\kap,a\} \subset x=\mathcal{H}_{a}(\kappa)
\to
\exi! b\leq\kap(b=\psi^{\vec{\xi}}_{\kap}(a))
)$, where
$lh(\vec{\xi})=N-2$.
\\

Since 
the cardinality of the set $\mathcal{H}_{\veps_{\mK+1}}(\pi)$ is $\pi$ for any infinite
cardinal $\pi\leq\mK$,
pick an injection $f :\mathcal{H}_{\Lam}(\mK)\to \mK$ so that
$f "\mathcal{H}_{\Lam}(\pi)\subset\pi$ for any weakly inaccessibles $\pi\leq\mK$.

\begin{lemma}\label{lem:welldefinedness}

\benu

\item\label{lem:welldefinedness.2}
$\forall a<\Lam\, A(a)
$.

\item\label{lem:welldefinedness.3}

$\pi\in Mh^{a}_{k}(\xi)$ is a $\Pi^{1}_{k-1}$-class on $L_{\pi}$ uniformly for weakly inaccessible cardinals $\pi\leq\mK$
and $a,\xi$.
This means that for each $k$
 there exists a $\Pi^{1}_{k-1}$-formula $mh^{a}_{k}(x)$ such that
$\pi\in Mh^{a}_{k}(\xi)$ iff $L_{\pi}\models mh^{a}_{k}(\xi)$ 
for any weakly inaccessible cardinals $\pi\leq\mK$ with $f"(\{a\}\cup K(\xi))\subset L_{\pi}$.

\item\label{lem:welldefinedness.1}
$\mK\in Mh^{\alp}_{N-1}(\Lam)\cap M_{N-1}(Mh^{\alp}_{N-1}(\Lam))$.

\eenu
\end{lemma}
\bprf
\\
\ref{lem:welldefinedness}.\ref{lem:welldefinedness.2}.
We show that $A(a)$ is progressive, i.e.,
$\forall a<\Lam[\forall c<a\, A(c) \to A(a)]$.

Assume $\forall c<a\, A(c)$ and $a<\Lam$.
$\forall b<\mK \exists ! x[x=\mathcal{H}_{a}(b)]$
 follows from IH in ${\sf ZFL}$.
$\exists! b\leq\kap(b=\psi^{\vec{\xi}}_{\kappa}a)$ follows from this.
\\

\noindent
\ref{lem:welldefinedness}.\ref{lem:welldefinedness.3}.
Let $\pi$ be a weakly inaccessible cardinal with $f"(\{a\}\cup K(\xi))\subset L_{\pi}$.
Let $f$ be an injection such that $f"\calh_{\Lam}(\pi)\subset L_{\pi}$.
Then for $\fal\alp\in K(\xi)(f(\alp)\in f" \calh_{\alp}(\pi))$,
$\pi\in Mh^{a}_{k}(\xi)$ iff
for any $f(\vec{\nu})=(f(\nu_{k}),\ldots,f(\nu_{N-1}))$, each of $f(\nu_{i})\in L_{\pi}$,
if $\fal \alp\in K(\vec{\nu})(f(\alp)\in f"\calh_{a}(\pi))$ and 
$\vec{\nu}<\xi$, then $\pi\in M_{k}(Mh^{a}_{k}(\vec{\nu}))$, where
$f"\calh_{a}(\pi)\subset L_{\pi}$ is a class in $L_{\pi}$.
\\

\noindent
\ref{lem:welldefinedness}.\ref{lem:welldefinedness.1}.
We show the following $B(a)$ is progressive in $a<\Lam$:
\[
B(a)  :\Lrarw 
 \mK\in Mh^{\alp}_{N-1}(a)\cap M_{N-1}(Mh^{\alp}_{N-1}(a))
\]
Note that $a\in \mathcal{H}_{a}(\mK)$ holds for any $a<\Lam$.

Suppose $\forall b<a\, B(b)$.
We have to show that $Mh^{\alp}_{N-1}(a)$ is $\Pi^{1}_{N-3}$-indescribable in $\mK$.
It is easy to see that if
$\pi\in M_{N-1}(Mh^{\alp}_{N-1}(a))$, then $\pi\in Mh^{\alp}_{N-1}(a)$
by induction on $\pi$.
Let $\tht(u)$ be a $\Pi^{1}_{N-3}$-formula such that $L_{\mK}\models\tht(u)$.

By IH we have $\forall b<a[\mK\in M_{N-1}(Mh^{\alp}_{N-1}(b))]$.
In other words, $\mK\in Mh^{\alp}_{N-1}(a)$, i.e., $L_{\mK}\models mh^{\alp}_{N-1}(a)$, where
$mh^{\alp}_{N-1}(a)$ is a $\Pi^{1}_{N-2}$-sentence in Proposition \ref{lem:welldefinedness}.\ref{lem:welldefinedness.3}.
Since the universe $L_{\mK}$ is $\Pi^{1}_{N-2}$-indescribable, pick a $\pi<\mK$ such that
$L_{\pi}$ enjoys the $\Pi^{1}_{N-2}$-sentence $\tht(u)\land mh^{\alp}_{N-1}(a)$, 
and $\{f(\alp),f(a)\}\subset L_{\pi}$.
Therefore  $\pi\in Mh^{\alp}_{N-1}(a)$ and $L_{\pi}\models\tht(u)$.
Thus $\mK\in M_{N-1}(Mh^{\alp}_{N-1}(a))$.
\eprf

\subsection{Normal forms in ordinal notations}
In this subsection we introduce an \textit{irreducibility} of sequences,
which is needed to define a normal form in ordinal notations.



\bprp\label{prp:Mhelementary}
$\pi\in Mh^{a}_{k}(\zeta) \spand \xi\leq\zeta \Rarw \pi\in Mh^{a}_{k}(\xi)$.

\eprp
\bprf
 (\ref{eq:dfMhkh}) for $\pi\in Mh^{a}_{k}(\xi)$ in Definition \ref{df:Cpsiregularsm}.\ref{df:Cpsiregularsm.2}
 follows from $\pi\in Mh^{a}_{k}(\zeta)$ and 
 Proposition \ref{prp:idless}.
\eprf

\blem\label{lem:stepdown}{\rm (Cf. Lemma 3 in \cite{LMPS}.)}
Assume $\mK\geq\pi\in Mh^{a}_{k}(\xi)\cap Mh^{a}_{k+1}(\xi_{0})$ with 
$2\leq k\leq N-1$, 
$he(\mu)\leq\xi_{0}$ and
$\{a\}\cup K(\mu)\subset\calh_{a}(\pi)$.
Then
$\pi\in Mh^{a}_{k}(\xi+\mu)$ holds.
Moreover if $\pi\in M_{k+1}$, then
$\pi\in M_{k+1}(Mh^{a}_{k}(\xi+\mu))$
holds.
\elem
\bprf
Suppose $\pi\in Mh^{a}_{k}(\xi)\cap Mh^{a}_{k+1}(\xi_{0})$ and $K(\mu)\subset\calh_{a}(\pi)$ with $he(\mu)\leq\xi_{0}$.
We show $\pi\in Mh^{a}_{k}(\xi+\mu)$ 
by induction on ordinals $\mu$.
First note that if $b\in\calh_{a}(\pi)$, then $f(b)\in f"\calh_{\Lam}(\pi)\subset L_{\pi}$.
We have $K(\xi+\mu)\subset\calh_{a}(\pi)$.
$\pi\in M_{k+1}(Mh^{a}_{k}(\xi+\mu))$ follows from $\pi\in Mh^{a}_{k}(\xi+\mu)$ and $\pi\in M_{k+1}$.

Let 
$(\zeta)*\vec{\nu}<\xi+\mu$ 
and $K(\zeta)\cup K(\vec{\nu})\subset\calh_{a}(\pi)$ 
for $\vec{\nu}=(\nu_{0},\ldots,\nu_{n-1})$.
We need to show that $\pi\in M_{k}(Mh^{a}_{k}((\zeta)*\vec{\nu}))$.
By Definition \ref{df:Lam}.\ref{df:Exp2.5}, let $(\zeta_{0})*(\mu_{0},\ldots,\mu_{n-1})$ be a sequence such that
$\zeta<\zeta_{0}\leq_{pt}\xi+\mu$, $\mu_{0}\leq_{pt}te(\zeta_{0})$, $\fal i\leq n-1(\nu_{i}<\mu_{i})$, 
and $\fal i<n-1(\mu_{i+1}\leq_{pt} te(\mu_{i}))$.

If $\zeta_{0}\leq_{pt}\xi$, then $(\zeta)*\vec{\nu}<\xi$, and $\pi\in M_{k}(Mh^{a}_{k}((\zeta)*\vec{\nu}))$
by $\pi\in Mh^{a}_{k}(\xi)$.

Let $\zeta_{0}=\xi+\zeta_{1}$ with $0<\zeta_{1}\leq_{pt}\mu$.
If $\zeta_{1}<_{pt}\mu$, then by IH with $he(\zeta_{1})=he(\mu)$ we have
$\pi\in Mh^{a}_{k}(\zeta_{0})$.
On the other hand we have $(\zeta)*\vec{\nu}<\zeta_{0}$.
Hence $\pi\in M_{k}(Mh^{a}_{k}((\zeta)*\vec{\nu}))$.

Finally consider the case when $0<\zeta_{1}=\mu$.
Then we obtain $\vec{\nu}<te(\xi+\mu)=te(\mu)\leq he(\mu)\leq\xi_{0}$.
$\pi\in Mh^{a}_{k+1}(\xi_{0})$ with Proposition \ref{prp:Mhelementary} yields $\pi\in M_{k+1}(Mh^{a}_{k+1}(\vec{\nu}))$.

On the other side we see $\pi\in Mh^{a}_{k}(\zeta)$ as follows.
We have $\zeta<\xi+\mu$.
If $\zeta\leq\xi$, then this follows from $\pi\in Mh^{a}_{k}(\xi)$ and Proposition \ref{prp:Mhelementary},
and if $\zeta=\xi+\lam<\xi+\mu$, then IH yields $\pi\in Mh^{a}_{k}(\zeta)$.


Since $\pi\in Mh^{a}_{k}(\zeta)$ is a $\Pi^{1}_{k-1}$-sentence holding on $L_{\pi}$ 
by Lemma \ref{lem:welldefinedness}.\ref{lem:welldefinedness.3} and $\{a\}\cup K(\zeta)\subset\calh_{a}(\pi)$,
we obtain $\pi\in M_{k+1}(Mh^{a}_{k}((\zeta)*\vec{\nu}))$, a fortiori  $\pi\in M_{k}(Mh^{a}_{k}((\zeta)*\vec{\nu}))$.
\eprf



\bdf\label{df:nfform2}
{\rm
For sequences of ordinals $\vec{\xi}=(\xi_{k},\ldots,\xi_{N-1})$ and $\vec{\nu}=(\nu_{k},\ldots,\nu_{N-1})$
and $2\leq k,m,n\leq N-1$,
 \[
 Mh^{a}_{m}(\vec{\nu})\prec_{k} Mh^{a}_{n}(\vec{\xi})
 :\Lrarw 
\fal\pi\in Mh^{a}_{n}(\vec{\xi})(\{a,\pi\}\cup K(\vec{\nu})\subset\calh_{a}(\pi) \Rarw \pi\in M_{k}(Mh^{a}_{m}(\vec{\nu})))
 .\]
}
\edf

\bcor\label{cor:stepdown}
Let $\vec{\nu}$ be a sequence defined from a sequence $\vec{\xi}$ as follows.
$\fal i<k(\nu_{i}=\xi_{i})$, $\fal i>k(\nu_{i}=0)$, and
$\nu_{k}=\xi_{k}+\Lam^{\xi_{k+1}}b$, where $2\leq k<N$, $b<\Lam$ and $\xi_{k+1}\neq 0$.
Then
$Mh^{a}_{2}(\vec{\nu})\prec_{k+1}Mh^{a}_{2}(\vec{\xi})$ holds.
In particular if $\pi\in Mh^{a}_{2}(\vec{\xi})$
 and $K(\vec{\nu})\cup\{\pi,a\}\subset\calh_{a}(\pi)$, then
$\psi_{\pi}^{\vec{\nu}}(a)<\pi$.
\ecor
\bprf
This is seen from Lemma \ref{lem:stepdown}.
\eprf

\bprp\label{th:KM1}
Let $\vec{\nu}=(\nu_{2},\ldots,\nu_{N-1})$, $\vec{\xi}=(\xi_{2},\ldots,\xi_{N-1})$ be sequences of ordinals$<\veps(\Lam)$ such that $\vec{\nu}<_{k}\vec{\xi}$ for
an integer $k$ with $2\leq k\leq N-1$.
Then $Mh^{a}_{2}(\vec{\nu})\prec_{k}Mh^{a}_{2}(\vec{\xi})$.
In particular if $\pi\in Mh^{a}_{2}(\vec{\xi})$
 and $K(\vec{\nu})\cup\{\pi,a\}\subset\calh_{a}(\pi)$, then
$\psi_{\pi}^{\vec{\nu}}(a)<\pi$.
\eprp
\bprf
Assume $\pi\in Mh^{a}_{2}(\vec{\xi})$ and $K(\vec{\nu})\subset\calh_{a}(\pi)$.
We have $\pi\in Mh^{a}_{k}(\xi_{k})$.  
By the definition (\ref{eq:dfMhkh}) and $(\nu_{k},\ldots,\nu_{N-1})<\xi_{k}$,
we obtain 
$\pi\in M_{k}(\bigcap_{k\leq i\leq N-1}Mh^{a}_{i}(\nu_{i}))$.

On the other hand we have $\pi\in \bigcap_{i<k}Mh^{a}_{i}(\xi_{i})$, and hence 
$\pi\in \bigcap_{i<k}Mh^{a}_{i}(\nu_{i})$
by  $\fal i<k(\nu_{i}\leq\xi_{i})$ and Proposition \ref{prp:Mhelementary}.
Since $\pi\in \bigcap_{i<k}Mh^{a}_{i}(\nu_{i})$ is a $\Pi^{1}_{k-2}$-sentence holding in $L_{\pi}$,
we obtain $\pi\in M_{k}(\bigcap_{i\leq N-1}Mh^{a}_{i}(\nu_{i}))=M_{k}(Mh^{a}_{2}(\vec{\nu}))$, 
a fortiori $\pi\in M_{2}(Mh^{a}_{2}(\vec{\nu}))$.

Suppose $\{\pi,a\}\subset\calh_{a}(\pi)$.
The set
$C=\{\kap<\pi: \calh_{a}(\kap)\cap\pi\subset\kap, K(\vec{\nu})\cup\{\pi,a\}\subset\calh_{a}(\kap)\}$
 is a club subset of the regular cardinal $\pi$.
This shows the existence of a $\kap\in Mh_{2}^{a}(\vec{\nu})\cap C\cap\pi$, and hence
$\psi_{\pi}^{\vec{\nu}}(a)<\pi$ by the definition (\ref{eq:Psivec}).
\eprf

\bprp\label{prp:nfform}
Let $\vec{\xi}=(\xi_{2},\ldots,\xi_{N-1})$ be a sequence of ordinals$<\veps(\Lam)$ such that
$\{\pi,a\}\cup K(\vec{\xi})\subset\calh_{a}(\pi)$.
Assume $Tl(\xi_{i})<\Lam_{k}(\xi_{i+k}+1)$ for some $i<N-1$ and $k>0$.
Then
$\pi\in Mh^{a}_{2}(\vec{\xi}) \Lrarw \pi\in Mh^{a}_{2}(\vec{\mu})$,
where
$\vec{\mu}=(\mu_{2},\ldots,\mu_{N-1})$ with
$\mu_{i}=\xi_{i}-Tl(\xi_{i})$ and $\mu_{j}=\xi_{j}$ for $j\neq i$.
\eprp
\bprf
When $0<\xi_{i}=\Lam^{\gam_{m}}a_{m}+\cdots+\Lam^{\gam_{1}}a_{1}+\Lam^{\gam_{0}}a_{0}$ with
$\gam_{m}>\cdots>\gam_{1}>\gam_{0}, \, 0<a_{i}<\Lam$,
$\mu_{i}=\Lam^{\gam_{m}}a_{m}+\cdots+\Lam^{\gam_{1}}a_{1}$ for $Tl(\xi_{i})=\Lam^{\gam_{0}}a_{0}$.
If $\xi_{i}=0$, then so is $\mu_{i}=0$.

Let $\pi\in Mh^{a}_{2}(\vec{\mu})$ and $Tl(\xi_{i})<\Lam_{k}(\xi_{i+k}+1)$.
We obtain $\fal j\leq k(he^{(j)}(Tl(\xi_{i}))<\Lam_{k-j}(\xi_{i+k}+1))$, and $he^{(k)}(Tl(\xi_{i}))\leq\xi_{i+k}$.
On the other hand we have $\pi\in Mh^{a}_{i+k}(\xi_{i+k})$.
From Lemma \ref{lem:stepdown} we see inductively that for any $j< k$, $\pi\in Mh^{a}_{i+j}(he^{(j)}(Tl(\xi_{i})))$.
In particular $\pi\in Mh^{a}_{i+1}(he(Tl(\xi_{i})))$, and once again by Lemma \ref{lem:stepdown} and $\pi\in Mh^{a}_{i}(\mu_{i})$
we obtain $\pi\in Mh^{a}_{i}(\xi_{i})$.
Hence $\pi\in Mh^{a}_{2}(\vec{\xi})$.
\eprf

\bdf\label{df:nfform1}
{\rm
A sequence of ordinals $\vec{\xi}=(\xi_{2},\ldots,\xi_{N-1})$ is said to be \textit{irreducible}
iff
$\fal i<N-1\fal k>0(\xi_{i}>0 \Rarw Tl(\xi_{i})\geq\Lam_{k}(\xi_{i+k}+1))$.
 }
 \edf

\bprp\label{prp:irreducible}
Let $\vec{\nu}=(\nu_{k},\ldots,\nu_{N-1})\neq\vec{0}$ be an irreducible sequence,
and $k_{0}\geq k$ be the least number such that $\nu_{k_{0}}\neq 0$.
Assume $\nu_{k_{0}}<he^{(k_{0}-k)}(\xi)$. Then
$\vec{\nu}<\xi$ holds in the sense of Definition \ref{df:Lam}.\ref{df:Exp2.5}.
\eprp
\bprf
Let
$\ell< N-k$ be the largest number such that
$\nu_{k+\ell}\neq 0$.
We show $(\nu_{k},\ldots,\nu_{k+\ell})<\xi$.
Since $\vec{\nu}$ is irreducible, we have $\Lam_{i}(\nu_{k_{0}+i}+1)\leq Tl(\nu_{k_{0}})$.
From $\nu_{k_{0}}<he^{(k_{0}-k)}(\xi)$ and $te(\mu)\leq he(\mu)$ we obtain
 $\nu_{k_{0}+i}<\nu_{k_{0}+i}+1\leq he^{(i)}(\nu_{k_{0}})\leq he^{(k_{0}-k+i)}(\xi)$.
Let $(\mu_{k},\ldots,\mu_{N-1})\subset_{pt}\xi$ such that $\mu_{k}=Hd(\xi)$ and
$\mu_{i+1}=he(\mu_{i})=te(Hd(\mu_{i}))$.
Then $te(\mu_{k+i})=he(\mu_{k+i})$ and $\mu_{k_{0}+i}=he(\mu_{k_{0}+i-1})=he^{(k_{0}-k+i)}(\xi)$ for $k_{0}-k+i>0$.
Therefore $(\mu_{k},\ldots,\mu_{k+\ell})\subset_{pt}\xi$ witnesses 
$(\nu_{k},\ldots,\nu_{k+\ell})<\xi$.
\eprf

\bdf\label{df:lx}
 {\rm Let}  $\vec{\xi}=(\xi_{k},\ldots,\xi_{N-1})$, $\vec{\nu}=(\nu_{k},\ldots,\nu_{N-1})$ 
 {\rm and} $\vec{\nu}\neq\vec{\xi}$.
  {\rm Let} $i\geq k$ {\rm be the minimal number such that}
 $\nu_{i}\neq\xi_{i}$. 
 {\rm Suppose} $(\xi_{i},\ldots,\xi_{N-1})\neq\vec{0}${\rm , and let}
  $k_{1}\geq i$ {\rm be the minimal number such that} $\xi_{k_{1}}\neq 0$. {\rm Then}
  $\vec{\nu}<_{lx,k}\vec{\xi}$ {\rm iff one of the followings holds:}
  \benu
  \item
  $(\nu_{i},\ldots,\nu_{N-1})=\vec{0}$.
  
  \item
  {\rm In what follows assume}  $(\nu_{i},\ldots,\nu_{N-1})\neq\vec{0}${\rm , and let}
  $k_{0}\geq i$ {\rm be the minimal number such that} $\nu_{k_{0}}\neq 0\, (i=\min\{k_{0},k_{1}\})$. {\rm Then}
  $\vec{\nu}<_{lx,k}\vec{\xi}$ {\rm iff one of the followings holds:}
  
   \benu
   \item
   $i=k_{0}< k_{1}$ {\rm and} $he^{(k_{1}-k_{0})}(\nu_{k_{0}})\leq\xi_{k_{1}}$.
   \item
   $k_{0}\geq k_{1}=i$ {\rm and} $\nu_{k_{0}}<he^{(k_{0}-k_{1})}(\xi_{k_{1}})$.
   \eenu
  \eenu

\edf

\bprp\label{prp:psinucomparison}
Suppose that both of $\vec{\nu}$ and $\vec{\xi}$ are irreducible. Then
$\vec{\nu}<_{lx,k}\vec{\xi} \Rarw
Mh^{a}_{k}(\vec{\nu})\prec_{k} Mh^{a}_{k}(\vec{\xi})
$.
\eprp
\bprf
Let $\pi\in Mh^{a}_{k}(\vec{\xi})$, $K(\vec{\nu})\subset\calh_{a}(\pi)$, 
and $i\geq k$ be the minimal number such that $\nu_{i}\neq\xi_{i}$.
We have $\pi\in\bigcap_{k\leq j<i} Mh^{a}_{j}(\nu_{j})$, which is a $\Pi^{1}_{i-2}$-sentence holding on $L_{\pi}$.
In the case $\xi_{i}\neq 0$, it suffices to show that $\pi\in M_{i}(\bigcap_{j\geq i} Mh^{a}_{j}(\nu_{j}))$,
 since then we obtain
$\pi\in M_{i}(Mh^{a}_{k}(\vec{\nu}))$ by $\pi\in Mh^{a}_{i}(\xi_{i})\subset M_{i}$, a fortiori $\pi\in M_{k}(Mh^{a}_{k}(\vec{\nu}))$.

If $(\nu_{i},\ldots,\nu_{N-1})=\vec{0}$, then $\xi_{i}\neq 0$ and
$\bigcap_{j\geq i} Mh^{a}_{j}(\nu_{j})$ denotes the class of limit ordinals.
Obviously $\pi\in M_{i}(\bigcap_{j\geq i} Mh^{a}_{j}(\nu_{j}))$.

In what follows assume  $(\nu_{i},\ldots,\nu_{N-1})\neq\vec{0}$, and let
$k_{0}\geq i$ be the minimal number such that $\nu_{k_{0}}\neq 0$,
and $k_{1}\geq i$ be the minimal number such that $\xi_{k_{1}}\neq 0$.
\\
{\bf Case 1}. $k_{0}\geq k_{1}=i$: Then we have $\nu_{k_{0}}<he^{(k_{0}-k_{1})}(\xi_{k_{1}})$.
Proposition \ref{prp:irreducible} yields
$(\nu_{k_{0}},\ldots,\nu_{N-1})<\xi_{k_{1}}=\xi_{i}$,
which in turn yields
 $\pi\in M_{i}(\bigcap_{j\geq i} Mh^{a}_{j}(\nu_{j}))$ by the definition (\ref{eq:dfMhkh}) of 
$\pi\in Mh^{a}_{i}(\xi_{i})$.
\\

\noindent
{\bf Case 2}. $i=k_{0}< k_{1}$:
Then we have $he^{(k_{1}-i)}(\nu_{i})\leq\xi_{k_{1}}$.
Also $\nu_{i+p}<he^{(p)}(\nu_{i})$ for any $p>0$ since $\vec{\nu}$ is irreducible and $\nu_{i}\neq 0$.
Let $j\geq k_{1}$.
Then $\nu_{j}<he^{(j-i)}(\nu_{i})\leq he^{(j-k_{1})}(\xi_{k_{1}})$.
In particular
$\nu_{k_{1}}<\xi_{k_{1}}$.
Proposition \ref{prp:irreducible} yields
$(\nu_{k_{1}},\ldots,\nu_{N-1})<\xi_{k_{1}}$.
$\pi\in Mh^{a}_{k_{1}}(\xi_{k_{1}})$ yields $\pi\in M_{k_{1}}(\bigcap_{j\geq k_{1}} Mh^{a}_{j}(\nu_{j}))$.
Moreover for any $p< k_{1}-i$, $he^{(k_{1}-i-p)}(\nu_{i+p})\leq \xi_{k_{1}}$ by Proposition \ref{prp:headcomparison}.
Lemma \ref{lem:stepdown} yields $\pi\in \bigcap_{k_{1}>j\geq i} Mh^{a}_{j}(\nu_{j})$.
Therefore $\pi\in M_{k_{1}}(Mh^{a}_{k}(\vec{\nu}))$, a fortiori $\pi\in M_{k}(Mh^{a}_{k}(\vec{\nu}))$.
\eprf

\bprp\label{prp:psicomparison}{\rm (Cf. Proposition 4.20 in \cite{Rathjen94})}\\
Let $\vec{\nu}=(\nu_{2},\ldots,\nu_{N-1})$, $\vec{\xi}=(\xi_{2},\ldots,\xi_{N-1})$ be 
irreducible sequences of ordinals$<\veps(\Lam)$, and assume that
$\psi_{\pi}^{\vec{\nu}}(b)<\pi$ and $\psi_{\kap}^{\vec{\xi}}(a)<\kap$.

Then $\bet_{1}=\psi_{\pi}^{\vec{\nu}}(b)<\psi_{\kap}^{\vec{\xi}}(a)=\alp_{1}$ iff one of the following cases holds:
\benu
\item\label{prp:psicomparison.0}
$\pi\leq \psi_{\kap}^{\vec{\xi}}(a)$.

\item\label{prp:psicomparison.1}
$b<a$, $\psi_{\pi}^{\vec{\nu}}(b)<\kap$ and 
$K(\vec{\nu})\cup\{\pi,b\}\subset\calh_{a}(\psi_{\kap}^{\vec{\xi}}(a))$.

\item\label{prp:psicomparison.2}
$b>a$ and $K(\vec{\xi})\cup\{\kap,a\}\not\subset\calh_{b}(\psi_{\pi}^{\vec{\nu}}(b))$.

\item\label{prp:psicomparison.25}
$b=a$, $\kap<\pi$ and $\kap\not\in\calh_{b}(\psi_{\pi}^{\vec{\nu}}(b))$.

\item\label{prp:psicomparison.3}
$b=a$, $\pi=\kap$, $K(\vec{\nu})\subset\calh_{a}(\psi_{\kap}^{\vec{\xi}}(a))$, and
$\vec{\nu}<_{lx,2}\vec{\xi}$.

\item\label{prp:psicomparison.4}
$b=a$, $\pi=\kap$, 
$K(\vec{\xi})\not\subset\calh_{b}(\psi_{\pi}^{\vec{\nu}}(b))$.

\eenu

\eprp
\bprf
If the case (\ref{prp:psicomparison.1}) holds, then 
$\psi_{\pi}^{\vec{\nu}}(b)\in\calh_{a}(\psi_{\kap}^{\vec{\xi}}(a))\cap\kap\subset\psi_{\kap}^{\vec{\xi}}(a)$.

If one of the cases (\ref{prp:psicomparison.2}) and (\ref{prp:psicomparison.25}) holds, then 
$K(\vec{\xi})\cup\{\kap,a\}\not\subset\calh_{a}(\psi_{\pi}^{\vec{\nu}}(b))$.
On the other hand we have 
$K(\vec{\xi})\cup\{\kap,a\}\subset\calh_{a}(\psi_{\kap}^{\vec{\xi}}(a))$.
Hence $\psi_{\pi}^{\vec{\nu}}(b)<\psi_{\kap}^{\vec{\xi}}(a)$.

If the case (\ref{prp:psicomparison.3}) holds, then Proposition \ref{prp:psinucomparison} yields
$Mh^{a}_{2}(\vec{\nu})\prec_{2} Mh^{a}_{2}(\vec{\xi})\ni\psi_{\kap}^{\vec{\xi}}(a)$.
Hence $\psi_{\kap}^{\vec{\xi}}(a)\in M_{2}(Mh^{a}_{2}(\vec{\nu}))$.
Since 
$K(\vec{\nu})\cup\{\kap,a\}\subset\calh_{a}(\psi_{\kap}^{\vec{\xi}}(a))$,
the set $\{\rho<\psi_{\kap}^{\vec{\xi}}(a): \calh_{a}(\rho)\cap\kap\subset\rho, K(\vec{\nu})\cup\{\kap,a\}\subset\calh_{a}(\rho)\}$ is club in $\psi_{\kap}^{\vec{\xi}}(a)$.
Therefore
$\psi_{\pi}^{\vec{\nu}}(b)=\psi_{\kap}^{\vec{\nu}}(a)<\psi_{\kap}^{\vec{\xi}}(a)$
by (\ref{eq:Psivec}) in Definition \ref{df:Cpsiregularsm}.\ref{df:Cpsiregularsm.3}.

Finally assume that the case (\ref{prp:psicomparison.4}) holds.
Since $K(\vec{\xi})\subset\calh_{a}(\psi_{\kap}^{\vec{\xi}}(a))$,
$\psi_{\pi}^{\vec{\nu}}(b)<\psi_{\kap}^{\vec{\xi}}(a)$ holds.

Conversely assume that $\psi_{\pi}^{\vec{\nu}}(b)<\psi_{\kap}^{\vec{\xi}}(a)$ and
$\psi_{\kap}^{\vec{\xi}}(a)<\pi$.

First consider the case $b<a$.
Then we have $K(\vec{\nu})\cup\{\pi,b\}\subset\calh_{b}(\psi_{\pi}^{\vec{\nu}}(b))\subset
\calh_{a}(\psi_{\kap}^{\vec{\xi}}(a))$. Hence (\ref{prp:psicomparison.1}) holds.

Next consider the case $b>a$.
$K(\vec{\xi})\cup\{\kap,a\}\subset\calh_{b}(\psi_{\pi}^{\vec{\nu}}(b))$ would yield
$\psi_{\kap}^{\vec{\xi}}(a)\in\calh_{b}(\psi_{\pi}^{\vec{\nu}}(b))\cap\pi\subset\psi_{\pi}^{\vec{\nu}}(b)$, a contradiction $\psi_{\kap}^{\vec{\xi}}(a)<\psi_{\pi}^{\vec{\nu}}(b)$.
Hence (\ref{prp:psicomparison.2}) holds.

Finally assume $b=a$.
Consider the case $\kap<\pi$.
$\kap\in \calh_{b}(\psi_{\pi}^{\vec{\nu}}(b))\cap\pi$ would yield 
$\psi_{\kap}^{\vec{\xi}}(a)<\kap<\psi_{\pi}^{\vec{\nu}}(b)$, a contradiction.
Hence $\kap\not\in \calh_{b}(\psi_{\pi}^{\vec{\nu}}(b))$, and 
(\ref{prp:psicomparison.25}) holds.
If $\pi<\kap$, then $\pi\in \calh_{b}(\psi_{\pi}^{\vec{\nu}}(b))\cap\kap
\subset\calh_{a}(\psi_{\kap}^{\vec{\xi}}(a))\cap\kap$, and $\pi<\psi_{\kap}^{\vec{\xi}}(a)$,
a contradiction, or we should say that (\ref{prp:psicomparison.0}) holds.
Finally let $\pi=\kap$.
We can assume that
$K(\vec{\xi})\subset\calh_{b}(\psi_{\pi}^{\vec{\nu}}(b))$, otherwise 
(\ref{prp:psicomparison.4}) holds.
If $\vec{\xi}<_{lx,2}\vec{\nu}$, then by (\ref{prp:psicomparison.3}) 
$\psi_{\kap}^{\vec{\xi}}(a)<\psi_{\pi}^{\vec{\nu}}(b)$ would follow.
If $K(\vec{\nu})\not\subset\calh_{a}(\psi_{\kap}^{\vec{\xi}}(a))$, then by 
(\ref{prp:psicomparison.4}) 
again $\psi_{\kap}^{\vec{\xi}}(a)<\psi_{\pi}^{\vec{\nu}}(b)$ would follow.
Hence $K(\vec{\nu})\subset\calh_{a}(\psi_{\kap}^{\vec{\xi}}(a))$
and $\vec{\nu}\leq_{lx}\vec{\xi}$.
If $\vec{\nu}=\vec{\xi}$, then $\psi_{\kap}^{\vec{\xi}}(a)=\psi_{\pi}^{\vec{\nu}}(b)$.
Therefore (\ref{prp:psicomparison.3}) must be the case.
\eprf
\\

Definition \ref{df:SD} is utilized to define a computable notation system in the next section \ref{subsec:decidable}.

\bdf\label{df:SD}
{\rm
A set $SD$ of sequences $\vec{\xi}=(\xi_{2},\ldots,\xi_{N-1})$ of ordinals $\xi_{i}<\veps(\Lam)$
is defined recursively as follows.
\benu
\item
$\vec{0}*(a)\in SD$ for each $a<\Lam$.
\item
(Cf.\,Definition \ref{df:Lam}.\ref{df:Exp2.10}.)
Let $\vec{\xi}=(\xi_{2},\ldots,\xi_{N-1})\in SD$, $1\leq k<N-1$, $\zeta<\veps(\Lam)$ be an ordinal such that 
$(\xi_{k+1},\ldots,\xi_{N-1})<_{sd}\zeta$, and
$(\xi_{2},\ldots,\xi_{k-1},\xi_{k},\zeta)*\vec{0}\in SD$.
Then for
 $\zeta_{k}=\xi_{k}+\Lam^{\zeta}a$ with an ordinal $a<\Lam$,
$(\xi_{2},\ldots,\xi_{k-1})*(\zeta_{k})*(\xi_{k+1},\ldots,\xi_{N-1})\in SD$
 and $(\xi_{2},\ldots,\xi_{k-1})*(\zeta_{k})*\vec{0}\in SD$.
\eenu
}
\edf

\bprp\label{prp:SD}
Let $\vec{\xi}=(\xi_{2},\ldots,\xi_{N-1})\in SD$.
\benu
\item\label{prp:SD.-1}
$(\xi_{2},\ldots,\xi_{i})*\vec{0}\in SD$ for each $i$ with $1\leq i<N$.
\item\label{prp:SD.0}
For $2\leq i<j<k<N$, if $\xi_{i}\neq 0$ and $\xi_{k}\neq 0$, then $\xi_{j}\neq 0$.

\item\label{prp:SD.1}
Let $\xi_{i}\neq 0$.
Then 
$(\xi_{i+1},\ldots,\xi_{N-1})<_{sd}te(\xi_{i})$.
\item\label{prp:SD.2}
$\vec{\xi}$ is irreducible.
\eenu
\eprp
\bprf
Let $1\leq k<N-1$, $\zeta<\veps(\Lam)$ be an ordinal such that 
$(\xi_{k+1},\ldots,\xi_{N-1})<_{sd}\zeta$, and
$(\xi_{2},\ldots,\xi_{k-1},\xi_{k},\zeta)*\vec{0}\in SD$.
Also let
 $\zeta_{k}=\xi_{k}+\Lam^{\zeta}a$ with an ordinal $a<\Lam$.
 \\
\ref{prp:SD}.\ref{prp:SD.-1} is seen by induction on the recursive definition of $\vec{\xi}\in SD$.
\\
\ref{prp:SD}.\ref{prp:SD.0} is seen by induction on the recursive definition of $\vec{\xi}\in SD$.
Suppose $\xi_{i}\neq 0$ for an $i<k$. 
From $(\xi_{2},\ldots,\xi_{k-1},\xi_{k},\zeta)*\vec{0}\in SD$ and $\zeta\neq 0$, IH yields $\xi_{k}\neq 0$.
\\
\ref{prp:SD}.\ref{prp:SD.1} and \ref{prp:SD}.\ref{prp:SD.2}.
We show these by simultaneous induction on the recursive definition of $\vec{\xi}\in SD$.
\\
\ref{prp:SD}.\ref{prp:SD.1}.
We show Proposition \ref{prp:SD}.\ref{prp:SD.1} for
the sequence $(\xi_{2},\ldots,\xi_{k-1})*(\zeta_{k})*(\xi_{k+1},\ldots,\xi_{N-1})\in SD$.
The proposition holds for the sequence $\vec{\xi}$, and we can assume $a\neq 0$.
We obtain $(\xi_{i+1},\ldots,\xi_{N-1})<_{sd}te(\xi_{i})$ for $i>k$ if $\xi_{i}\neq 0$, and
$(\xi_{k+1},\ldots,\xi_{N-1})<_{sd}te(\zeta_{k})=\zeta$ by the assumption.
Let $2\leq i<k$ and $\xi_{i}\neq 0$.
We show
$(\xi_{i+1},\ldots,\xi_{k-1})*(\zeta_{k})*(\xi_{k+1},\ldots,\xi_{N-1})<_{sd}te(\xi_{i})$.
It suffices to show that $\zeta_{k}<_{sd}te^{(k-i)}(\xi_{i})$.
By IH we have $\xi_{k}<_{sd}te^{(k-i)}(\xi_{i})$.
On the other hand we have $\xi_{k}\neq 0$ by $(\xi_{2},\ldots,\xi_{k-1},\xi_{k},\zeta)*\vec{0}\in SD$,
$\zeta\neq 0$, and Proposition \ref{prp:SD}.\ref{prp:SD.0}.
Moreover $(\xi_{2},\ldots,\xi_{k-1},\xi_{k},\zeta)*\vec{0}$ is irreducible by Proposition \ref{prp:SD}.\ref{prp:SD.2},
and hence $Tl(\xi_{k})\geq\Lam^{\zeta+1}$.
Therefore $te(\xi_{k})>\zeta$.
This means that $\zeta_{k}=_{NF}\xi_{k}+\Lam^{\zeta}a$, and $\xi_{k}<_{sd}te^{(k-i)}(\xi_{i})$
yields $\zeta_{k}<_{sd}te^{(k-i)}(\xi_{i})$ by Definition \ref{df:Lam}.\ref{df:Exp2.9}.
\\
\ref{prp:SD}.\ref{prp:SD.2}.
If $(\xi_{i+1},\ldots,\xi_{N-1})<_{sd}te(\xi_{i})$ for $\xi_{i}\neq 0$,
then $\xi_{i+k}<_{sd}te^{(k)}(\xi_{i})$ for $k>0$, and
$\xi_{i+k}+1\leq te^{(k)}(\xi_{i})$.
Hence $\Lam_{k}(\xi_{i+k}+1)\leq\Lam^{te(\xi_{i})}\leq Tl(\xi_{i})$, and
$\vec{\xi}$ is irreducible.
\eprf

\section{Computable notation system $OT$}\label{subsec:decidable}

In this section (except Propositions \ref{prp:l.5.4})
we work in a weak fragment of arithmetic, e.g., in the fragment $I\Sig_{1}$ or even in the bounded arithmetic $S^{1}_{2}$.
Referring Proposition \ref{prp:psicomparison}
the sets of ordinal terms $OT\subset\Lam=\veps_{\mK+1}$ 
and $E\subset\veps(\Lam)=\veps_{\mK+2}$
over symbols $\{0,\mK,\Lam,+,\ome,\vphi, \Ome, \psi\}$ are defined recursively.
$OT$ is isomorphic to a subset of $\calh_{\Lam}(0)$.
Simultaneously we define finite sets $K_{\del}(\alp)\subset OT$ for $\del,\alp\in OT$, and sequences
$(m_{k}(\alp))_{2\leq k\leq N-1}$ for $\alp\in OT\cap\mK$, where
in $\alp=\psi_{\pi}^{\vec{\nu}}(a)$, $m_{k}(\alp)=\nu_{k}$, i.e.,
$\vec{\nu}=(\nu_{2},\ldots,\nu_{N-1})=(m_{2}(\alp),\ldots,m_{N-1}(\alp))=(m_{k}(\alp))_{k}=\vec{m}(\alp)$.
For $\{\alp_{0},\ldots,\alp_{m},\bet\}\subset OT$,
$K_{\del}(\alp_{0},\ldots,\alp_{m}):=\bigcup_{i\leq m}K_{\del}(\alp_{i})$,
$K_{\del}(\alp_{0},\ldots,\alp_{m})<\bet:\Lrarw\fal\gam\in K_{\del}(\alp_{0},\ldots,\alp_{m})(\gam<\bet)$,
and 
$\bet\leq K_{\del}(\alp_{0},\ldots,\alp_{m}):\Lrarw \exi\gam\in K_{\del}(\alp_{0},\ldots,\alp_{m})(\bet\leq\gam)$.

An ordinal term in $OT$ is said to be a \textit{regular} term if it is one of the form $\mK$,
$\Ome_{\bet+1}$ or $\psi_{\pi}^{\vec{\nu}}(a)$ with the non-zero sequences 
$\vec{\nu}\neq\vec{0}$.
$\mK$ and the latter terms $\psi_{\pi}^{\vec{\nu}}(a)$ are \textit{Mahlo} terms.

$\alp=_{NF}\alp_{m}+\cdots+\alp_{0}$ means that $\alp=\alp_{m}+\cdots+\alp_{0}$ 
and $\alp_{m}\geq\cdots\geq\alp_{0}$
and each $\alp_{i}$ is a non-zero additive principal number.
$\alp=_{NF}\vphi\bet\gam$ means that $\alp=\vphi\bet\gam$ and $\bet,\gam<\alp$.
$\alp=_{NF}\ome^{\bet}$ means that $\alp=\ome^{\bet}>\bet$.
$\alp=_{NF}\Ome_{\bet}$ means that $\alp=\Ome_{\bet}>\bet$.


Let $pd(\psi_{\pi}^{\vec{\nu}}(a))=\pi$ (even if $\vec{\nu}=\vec{0}$).
Moreover for $n$,
$pd^{(n)}(\alp)$ is defined recursively by $pd^{(0)}(\alp)=\alp$ and
$pd^{(n+1)}(\alp)\simeq pd(pd^{(n)}(\alp))$.

For terms $\pi,\kap\in OT$,
$\pi\prec\kap$ denotes the transitive closure of the relation
$\{(\pi,\kap): \exi \vec{\xi}\exi b[\pi=\psi_{\kap}^{\vec{\xi}}(b)]\}$,
 and its reflexive closure
$\pi\preceq\kap:\Lrarw \pi\prec\kap \lor \pi=\kap\Lrarw\exi n(\kap=pd^{(n)}(\pi))$.

For each ordinal term $\alp=\psi_{\pi}^{\vec{\nu}}(a)$,
a series $(\pi_{i})_{i\leq L}$ of ordinal terms is uniquely determined as follows:
$\pi_{L}=\alp$, $\pi_{i}=pd(\pi_{i+1})$ and $\pi_{0}=\mK$.
Let us call the series $(\pi_{i})_{i\leq L}$ the \textit{collapsing series} of $\alp=\pi_{L}$.

Then we see that an ordinal term
$\alp=\psi_{\pi}^{\vec{\nu}}(a)$ with $\vec{\nu}\neq\vec{0}$ 
is constructed by Definition \ref{df:notationsystem}.\ref{df:notationsystem.10} below
iff $L=1$.
$\alp$ is constructed by Definition \ref{df:notationsystem}.\ref{df:notationsystem.12}
iff $L\equiv 1 \!\!\!\!\pmod{(N-2)}$.
Otherwise $\alp$ is constructed by Definition \ref{df:notationsystem}.\ref{df:notationsystem.11}.

\bdf\label{df:notationsystem}
{\rm
$\ell\alp$ denotes the number of occurrences of symbols
\\
$\{0,\mK,\Lam,+,\ome,\vphi, \Ome, \psi\}$
in terms $\alp\in OT\cup E$.

\benu
\item\label{df:notationsystem.1}

 \benu
 \item\label{df:notationsystem.3}
$0\in E$.

 \item\label{df:notationsystem.5-1}
If $0<a\in OT$,
then $a\in E$.
$K(a)=\{a\}$.

 \item\label{df:notationsystem.5}
If $\{\xi_{i}:i\leq m\}\subset E$, $\xi_{m}>\cdots>\xi_{0}>0$ and
$0<b_{i}\in OT$,
then $\sum_{i\leq m}\Lam^{\xi_{i}}b_{i}=
\Lam^{\xi_{m}}b_{m}+\cdots+\Lam^{\xi_{0}}b_{0}\in E$.
$K(\sum_{i\leq m}\Lam^{\xi_{i}}b_{i})=\{b_{i}:i\leq m\}\cup\bigcup\{K(\xi_{i}):i\leq m\}$.

\item
For sequences $\vec{\nu}=(\nu_{2},\ldots,\nu_{N-1})$, let
$K(\vec{\nu})=\bigcup_{2\leq i\leq N-1}K(\nu_{i})$.

  
 

 \eenu
 
 \item
 \benu
 \item\label{df:notationsystem.2}
$0,\mK\in OT$.
$m_{k}(0)=0$ for any $k$, and $K_{\del}(0)=K_{\del}(\mK)=\emptyset$.

 \item\label{df:notationsystem.4}
If $\alp=_{NF}\alp_{m}+\cdots+\alp_{0}\, (m>0)$ with $\{\alp_{i}:i\leq m\}\subset OT$, 
then
$\alp\in OT$, and $m_{k}(\alp)=0$ for any $k$.
$K_{\del}(\alp)=K_{\del}(\alp_{0},\ldots,\alp_{m})$.

 \item\label{df:notationsystem.6}
If $\alp=_{NF}\vphi\bet\gam$ with $\{\bet,\gam\}\subset OT\cap\mK$, then
$\alp\in OT$, and $m_{k}(\alp)=0$ for any $k$.
$K_{\del}(\alp)=K_{\del}(\bet,\gam)$.

 \item\label{df:notationsystem.7}
If $\alp=_{NF}\ome^{\bet}$ with $\mK<\bet\in OT$, then $\alp\in OT$,
and $m_{k}(\alp)=0$ for any $k$.
$K_{\del}(\alp)=K_{\del}(\bet)$.

 \item\label{df:notationsystem.8}
If $\alp=_{NF}\Ome_{\bet}$ with $\bet\in OT\cap\mK$, then
$\alp\in OT$.
$m_{2}(\alp)=1, m_{k}(\alp)=0$ for any $k>2$ if $\bet$ is a successor ordinal.
Otherwise
$m_{k}(\alp)=0$ for any $k$.
In each case
$K_{\del}(\alp)=K_{\del}(\bet)$.

 \item\label{df:notationsystem.9}
Let $\alp=\psi_{\pi}(a):=\psi_{\pi}^{\vec{0}}(a)$ where $\pi$ is a regular term , i.e.,
either $\pi=\mK$ or
$\vec{m}(\pi)\neq\vec{0}$,
and  $K_{\alp}(\pi,a)<a$.

Then
$\alp=\psi_{\pi}(a)\in OT$. Let
$m_{k}(\alp)=0$ for any $k$.
$K_{\del}(\psi_{\pi}(a))=\emptyset$ if $\alp<\del$.
$K_{\del}(\psi_{\pi}(a))=\{a\}\cup K_{\del}(a,\pi)$ otherwise.

 \item\label{df:notationsystem.10}
Let $\alp=\psi_{\mK}^{\vec{\nu}}(a)$ with
 $\vec{\nu}=\vec{0}*(b)\, (lh(\vec{\nu})=N-2)$ and $b,a\in OT$
such that $0<b\leq a$ and $K_{\alp}(b,a)<a$.

 Then
$\alp=\psi_{\mK}^{\vec{\nu}}(a)\in OT$.
Let
$m_{N-1}(\alp)=b$, $m_{k}(\alp)=0$
for $k<N-1$.
$K_{\del}(\psi_{\mK}^{\vec{\nu}}(a))=\emptyset$ if $\alp<\del$.
$K_{\del}(\psi_{\mK}^{\vec{\nu}}(a))=\{a\}\cup \bigcup\{K_{\del}(\gam):\gam\in K(\nu)\}$ otherwise.

 \item\label{df:notationsystem.11}
Let $\pi\in OT\cap\mK$ be such that 
$m_{k+1}(\pi)\neq 0$ 
and $\fal i>k+1(m_{i}(\pi)=0)$
for a $k\, (2\leq k\leq N-2)$, and
$b,a\in OT$ such that $0< b\leq a$.
Let $\vec{\nu}=(\nu_{2},\ldots,\nu_{N-1})$ be a sequence
defined by $\fal i<k(\nu_{i}=m_{i}(\pi))$,
$\nu_{k}=m_{k}(\pi)+\Lam^{m_{k+1}(\pi)}b$,
and $\fal i>k(\nu_{i}=0)$.

Then 
$\alp=\psi_{\pi}^{\vec{\nu}}(a)\in OT$ if
$K_{\alp}(\pi,a,b)\cup K_{\alp}(K(\vec{m}(\pi)))<a$.
Let $m_{i}(\alp)=\nu_{i}$ for each $i$.
$K_{\del}(\psi_{\pi}^{\vec{\nu}}(a))=\emptyset$ if $\alp<\del$.
Otherwise
$K_{\del}(\psi_{\pi}^{\vec{\nu}}(a))=
\{a\}\cup K_{\del}(a,\pi)\cup\bigcup\{K_{\del}(b): b\in K(\vec{\nu})\}$.

 \item\label{df:notationsystem.12}
Let $\pi\in OT\cap\mK$ be such that 
$m_{2}(\pi)\neq 0$ and $\fal i>2(m_{i}(\pi)=0)$,
and $a\in OT$.
Let $\vec{0}\neq\vec{\nu}=(\nu_{2},\ldots,\nu_{N-1})\in SD$ be a sequence
 of ordinal terms $\nu_{i}\in E$
such that $\vec{\nu}<_{sp} m_{2}(\pi)$.

Then 
$\alp=\psi_{\pi}^{\vec{\nu}}(a)$ if $K_{\alp}(\pi,a)<a$, and
\beqn\label{eq:notationsystem.12}
\forall k( K_{\alpha}(\nu_{k})<\max K(\nu_{k}))
\eeqn

Let $m_{i}(\alp)=\nu_{i}$ for each $i$.

$K_{\del}(\psi_{\pi}^{\vec{\nu}}(a))=\emptyset$ if $\alp<\del$.
Otherwise
$K_{\del}(\psi_{\pi}^{\vec{\nu}}(a))=\{a\}\cup K_{\del}(a,\pi)\cup\bigcup\{K_{\del}(b): b\in K(\vec{\nu})\}$.

 \eenu

\eenu
}
\edf

Let $\{\pi,a,\xi\}\subset \calh_{a}(\pi)$.
Then $\xi=m_{k}(\pi)$ is intended to be equivalent to $\pi\in Mh^{a}_{k}(\xi)$.
For Definition \ref{df:notationsystem}.\ref{df:notationsystem.11}, see Corollary \ref{cor:stepdown},
and for
Definition \ref{df:notationsystem}.\ref{df:notationsystem.12}, see Proposition \ref{th:KM1}.

\bprp\label{prp:SDm}
For each Mahlo term $\alp=\psi_{\pi}^{\vec{\nu}}(a)\in OT$,
$\vec{m}(\alp)=\vec{\nu}\in SD$ for the class $SD$ in
Definition \ref{df:SD}.
\eprp


\bprp\label{prp:l.5.4}
For any $\alp\in OT$ and any $\del$ such that $\del=0,\mK$ or $\del=\psi_{\pi}^{\vec{\nu}}(b)$ for some $\pi,b,\vec{\nu}$,
$\alp\in\calh_{\gam}(\del) \Lrarw K_{\del}(\alp)<\gam$.
\eprp
\bprf
By induction on $\ell\alp$.
\eprf


\blem\label{lem:compT}
$(OT,<)$ is a computable notation system of ordinals.
In particular the order type of the initial segment $\{\alp\in OT: \alp<\Ome_{1}\}$
is less than $\ome_{1}^{CK}$.

Specifically
each of $\alp<\bet$ and $\alp=\bet$ is decidable for $\alp,\bet\in OT$, 
and $\alp\in OT$ is decidable for terms $\alp$
over symbols $\{0,\mK,\Lam,+,\ome,\vphi, \Ome, \psi\}$.
\elem
\bprf
These are shown simultaneously referring Propositions \ref{prp:psicomparison} and \ref{prp:l.5.4}.
Let us give recursive definitions only for terms $\Ome_{\alp},\psi_{\kap}^{\vec{\nu}}(a)\in OT$.

First $\Ome_{\psi_{\kap}^{\vec{\nu}}(a)}=\psi_{\kap}^{\vec{\nu}}(a)$, i.e.,
$\Ome_{\alp}<\psi_{\kap}^{\vec{\nu}}(a) \Lrarw \alp<\psi_{\kap}^{\vec{\nu}}(a)$,
$\psi_{\kap}^{\vec{\nu}}(a)<\Ome_{\alp}\Lrarw \psi_{\kap}^{\vec{\nu}}(a)<\alp$.
Next $\Ome_{\alp}<\psi_{\Ome_{\alp+1}}(a)<\Ome_{\alp+1}$.

Finally for $\psi_{\pi}^{\vec{\nu}}(b),\psi_{\kap}^{\vec{\xi}}(a)\in OT$,
$\psi_{\pi}^{\vec{\nu}}(b)<\psi_{\kap}^{\vec{\xi}}(a)$ iff one of the following cases holds:
\benu
\item
$\pi\leq \psi_{\kap}^{\vec{\xi}}(a)$.

\item
$b<a$, $\psi_{\pi}^{\vec{\nu}}(b)<\kap$, and
$K_{\psi_{\kap}^{\vec{\xi}}(a)}(\{\pi,b\}\cup K(\vec{\nu}))<a$.

\item
$b\geq a$, and 
$b\leq K_{\psi_{\pi}^{\vec{\nu}}(b)}(\{\kap,a\}\cup K(\vec{\xi}))$.




\item
$b=a$, $\pi=\kap$, 
$K_{\psi_{\kap}^{\vec{\xi}}(a)}(K(\vec{\nu}))<a$, and
$\vec{\nu}<_{lx,2}\vec{\xi}$.

\eenu
\eprf



\bprp\label{prp:G6}
\benu
\item\label{prp:G6.3}
Let $\bet=\psi_{\pi}^{\vec{\nu}}(b)$ with $\pi=\psi_{\kap}^{\vec{\xi}}(a)$.
Then $a<b$.

\item\label{prp:G6.5}
For $\alpha=\psi_{\pi}^{\vec{\nu}}(a)\in OT$,
$\max K(\vec{\nu})\leq a$ holds.
\eenu
\eprp
\bprf
\ref{prp:G6}.\ref{prp:G6.3}.
Let $\bet=\psi_{\pi}^{\vec{\nu}}(b)$ with $\pi=\psi_{\kap}^{\vec{\xi}}(a)$.
Then $K_{\bet}(\{\pi,b\}\cup K(\vec{\nu}))<b$.
On the other hand we have $\bet<\pi$.
Hence $a\in K_{\bet}(\pi)<b$.
\\

\noindent
\ref{prp:G6}.\ref{prp:G6.5}.
This is seen by induction on $\ell\alpha$.
Ww have $c<a$ by Proposition \ref{prp:G6}.\ref{prp:G6.3} when $\pi=\psi_{\sig}^{\vec{\mu}}(c)$

When $\alpha$ is constructed by Definition \ref{df:notationsystem}.\ref{df:notationsystem.11},
$\nu_{k}=m_{k}(\pi)+\Lam^{m_{k+1}(\pi)}b$ holds for $b\leq a$.
By IH we have $\max K(\vec{m}(\pi))\leq c<a$ when $\pi=\psi_{\sig}^{\vec{\mu}}(c)$.

Suppose $\alpha$ is constructed by Definition \ref{df:notationsystem}.\ref{df:notationsystem.12}.
We obtain $\vec{\nu}<_{sp}m_{2}(\pi)$, and hence
$\max K(\vec{\nu})\leq \max K(m_{2}(\pi))\leq c<a$ by IH.
\eprf

\section{Operator controlled derivations}\label{sect:controlledOme}

In this section,
operator controlled derivations are defined, which are introduced by W. Buchholz\cite{Buchholz}.

In this and the next sections except otherwise stated
 $\alp,\bet,\gam,\ldots,a,b,c,d,\ldots$ range over ordinal terms in $OT\subset \calh_{\Lam}(0)$,
$\xi,\zeta,\nu,\mu,\iota,\ldots$ range over ordinal terms in $E$,
$\vec{\xi},\vec{\zeta},\vec{\nu},\vec{\mu},\vec{\iota},\ldots$ range over finite sequences over ordinal
terms in $E$,
and $\pi,\kap,\rho,\sig,\tau,\lam,\ldots$ range over regular ordinal terms
$\mK$, $\Ome_{\bet+1}$, $\psi_{\pi}^{\vec{\nu}}(a)$
with $\vec{\nu}\neq\vec{0}$.
$Reg$ denotes the set of regular ordinal terms.
We write $\alp\in\calh_{a}(\bet)$ for $K_{\bet}(\alp)<a$.

\subsection{Classes of sentences}\label{subsec:classformula}

Following Buchholz\cite{Buchholz} let us introduce a language for ramified set theory $RS$.

\bdf
$RS$\textit{-terms} {\rm and their} \textit{levels} {\rm are inductively defined.}
\benu
\item
{\rm For each $\alp\in OT\cap\mK$, $L_{\alp}$ is an $RS$-term of level $\alp$.}

\item
{\rm If $\phi(x,y_{1},\ldots,y_{n})$ is a set-theoretic formula in the language $\{\in\}$, and $a_{1},\ldots,a_{n}$ are
$RS$-terms of levels$<\alp$, then 
$[x\in L_{\alp}:\phi^{L_{\alp}}(x,a_{1},\ldots,a_{n})]$ 
is an $RS$-term of level $\alp$.}
\eenu

\edf
Each ordinal term $\alp$ is denoted by the ordinal term $[x\in L_{\alp}: x \mbox{ is an ordinal}]$, whose level is $\alp$.

\bdf\benu
\item
$|a|$ {\rm denotes the level of $RS$-terms $a$, and $Tm(\alp)$ the set of $RS$-terms of level$<\alp$.
$Tm=Tm(\mK)$ is then the set of $RS$-terms, which are denoted by $a,b,c,d,\ldots$}

\item
{\rm $RS$-\textit{formulas} 
are constructed from \textit{literals}
$a\in b, a\not\in b$
by propositional connectives $\lor,\land$, bounded quantifiers $\exi x\in a, \fal x\in a$ and
unbounded quantifiers $\exi x,\fal x$.
Unbounded quantifiers $\exists x,\forall x$ are denoted by $\exists x\in L_{\mK},\forall x\in L_{\mK}$, resp.}

\item
{\rm For $RS$-terms and $RS$-formulas $\iota$,
${\sf k}(\iota)$ denotes the set of ordinal terms $\alp$ such that the constant $L_{\alp}$ occurs in $\iota$.}

\item
{\rm For a set-theoretic $\Sig_{n}$-formula $\psi(x_{1},\ldots,x_{m})$ in $\{\in\}$ and $a_{1},\ldots,a_{m}\in Tm(\kap)$,
$\psi^{L_{\kap}}(a_{1},\ldots,a_{m})$ is a $\Sig_{n}(\kap)$\textit{-formula}, where
$n=0,1,2,\ldots$ and $\kap\leq\mK$. $\Pi_{n}(\kap)$-formulas are defined dually.}

\item
{\rm For} $\tht\equiv\psi^{L_{\kap}}(a_{1},\ldots,a_{m})\in\Sig_{n}(\kap)$ 
{\rm and} $\lam<\kap$,
$\tht^{(\lam,\kap)}:\equiv \psi^{L_{\lam}}(a_{1},\ldots,a_{m})$.
\eenu
\edf
Note that the level $|t|=\max(\{0\}\cup{\sf k}(t))$ for $RS$-terms $t$.
In what follows we need to consider \textit{sentences}.
Sentences are denoted $A,C$ possibly with indices.

The assignment of disjunctions and conjunctions to sentences is defined as in
\cite{Buchholz}.

\begin{definition}\label{df:assigndc}
\benu
\item
{\rm For} $b,a\in Tm(\mK)$ {\rm with} $|b|<|a|$,
\[
(b\veps a) :\equiv
\left\{
\begin{array}{ll}
A(b) & \mbox{{\rm if }} a\equiv[x\in L_{\alp}: A(x)]
\\
b\not\in L_{0} & \mbox{{\rm if }} a\equiv L_{\alp}
\end{array}
\right.
\]
{\rm and} 
$(a=b):\equiv(\fal x\in a(x\in b)\land \fal x\in b(x\in a))$.

\item\label{df:assigndc0}
{\rm For} $b,a\in Tm(\mK)$ {\rm and} $J:=Tm(|a|)$
\[
(b\in a):\simeq \bigvee(c\veps a \land c=b)_{c\in J}
\mbox{ {\rm and }}
(b\not\in a):\simeq \bigwedge(c\not\veps a \lor c\neq b)_{c\in J}
\]

\item
$(A_{0}\lor A_{1}):\simeq\bigvee(A_{\iota})_{\iota\in J}$
{\rm and}
$(A_{0}\land A_{1}):\simeq\bigwedge(A_{\iota})_{\iota\in J}$
{\rm for} $J:=2$.

\item
{\rm For} $a\in Tm(\mK)\cup\{L_{\mK}\}$ {\rm and} $J:=Tm(|a|)$
\[
\exists x\in a\, A(x):\simeq\bigvee(b\veps a \land A(b))_{b\in J}
\mbox{ {\rm and }}
\forall x\in a\, A(x):\simeq\bigwedge(b\not\veps a \lor A(b))_{b\in J}
.\]

\eenu

\end{definition}

The rank $\mbox{{\rm rk}}(\iota)$ of sentences or terms $\iota$ 
is defined as in \cite{Buchholz}.

\begin{definition}\label{df:rank}
\benu

\item\label{df:rank1}
$\mbox{{\rm rk}}(\lnot A):=\mbox{{\rm rk}}(A)$.

\item\label{df:rank2}
$\mbox{{\rm rk}}(L_{\alp})=\ome\alp$.

\item\label{df:rank3}
$\rk([x\in L_{\alp}: A(x)])=\max\{\ome\alp+1, \rk(A(L_{0}))+2\}$.

\item\label{df:rank4}
$\mbox{{\rm rk}}(a\in b)=\max\{\mbox{{\rm rk}}(a)+6,\mbox{{\rm rk}}(b)+1\}$.



\item\label{df:rank6}
$\mbox{{\rm rk}}(A_{0}\lor A_{1}):=\max\{\mbox{{\rm rk}}(A_{0}),\mbox{{\rm rk}}(A_{1})\}+1$.

\item\label{df:rank7}

 $\rk(\exists x\in a\, A(x)):=\max\{\omega \rk(a), \mbox{{\rm rk}}(A(L_{0}))+2\}$ {\rm for } $a\in Tm(\mK)\cup\{L_{\mK}\}$.

\eenu

\end{definition}

\begin{proposition}\label{lem:rank}
Let $A$ be a sentence with 
$A\simeq\bigvee(A_{\iota})_{\iota\in J}$ or $A\simeq\bigwedge(A_{\iota})_{\iota\in J}$.
\benu

\item\label{lem:rank0}
$\mbox{{\rm rk}}(A)<\mK+\omega$.

\item\label{lem:rank1}
$|A|\leq \mbox{{\rm rk}}(A)\in\{\omega |A|+i  : i\in\omega\}$.

\item\label{lem:rank2}
$\forall\iota\in J(\mbox{{\rm rk}}(A_{\iota})<\mbox{{\rm rk}}(A))$.

\item\label{lem:rank3}
$\rk(A)<\lam \Rarw A\in\Sig_{0}(\lam)$

\eenu

\end{proposition}

\subsection{Operator controlled derivations}\label{subsec:operatorcont}

{\rm By an} \textit{operator} {\rm we mean a map} $\mathcal{H}$, 
$\mathcal{H}:\mathcal{P}(OT)\to\mathcal{P}(OT)${\rm , such that}
\benu
\item
$\forall X\subset OT[X\subset\mathcal{H}(X)]$.

\item
$\forall X,Y\subset OT[Y\subset\mathcal{H}(X) \Rightarrow \mathcal{H}(Y)\subset\mathcal{H}(X)]$.
\eenu

{\rm For an operator} $\mathcal{H}$ {\rm and} $\Theta,\Tht_{1}\subset OT$,
$\mathcal{H}[\Theta](X):=\mathcal{H}(X\cup\Theta)${\rm , and}
$\mathcal{H}[\Theta][\Tht_{1}]:=(\mathcal{H}[\Theta])[\Tht_{1}]${\rm , i.e.,}
$\mathcal{H}[\Theta][\Tht_{1}](X)=\mathcal{H}(X\cup\Theta\cup\Tht_{1})$.

Obviously $\mathcal{H}_{\alpha}$ is an operator for any $\alp$, and
if $\calh$ is an operator, then so is $\mathcal{H}[\Theta]$.

\textit{Sequents} are finite sets of sentences, and inference rules are formulated in one-sided sequent calculus.
Let $\mathcal{H}=\calh_{\gam}\, (\gam\in OT)$ be an operator, $\Tht$ a finite set of $\mK$,
$\Gamma$ a sequent, 
$a\in OT$ and 
$b\in OT\cap(\mK+\omega)$.

We define a relation $(\mathcal{H}_{\gam},\Tht)\vdash^{a}_{b}\Gamma$, which is read `there exists an infinitary derivation
of $\Gamma$ which is $\Tht$\textit{-controlled} by $\mathcal{H}_{\gam}$, and 
 whose height is at most $a$ and its cut rank is less than $b$'.

$\kappa,\lambda,\sigma,\tau,\pi$ ranges over regular ordinal terms.

\begin{definition}\label{df:controlderreg}
{\rm
$(\mathcal{H}_{\gam},\Tht)\vdash^{a}_{b}\Gamma$ holds if

\begin{equation}\label{eq:controlder}
{\sf k}(\Gamma)\cup\{a\}\subset\mathcal{H}_{\gam}[\Tht] 
\end{equation}
and one of the following cases holds:

\bdes

\item[$(\bigvee)$]
$A\simeq\bigvee\{A_{\iota}: \iota\in J\}$, 
$A\in\Gamma$ and there exist $\iota\in J$
and
 $a(\iota)<a$ such that
\beqn\label{eq:bigveebnd}
|\iota|< a
\eeqn
and
$(\mathcal{H}_{\gam},\Tht)\vdash^{a(\iota)}_{b}\Gamma,A_{\iota}$.

\item[$(\bigwedge)$]
$A\simeq\bigwedge\{A_{\iota}: \iota\in J\}$, 
$A\in\Gamma$ and for every
$\iota\in J$ there exists an $a(\iota)<a$ 
such that
$(\mathcal{H}_{\gam},\Tht\cup\{{\sf k}(\iota)\})\vdash^{a(\iota)}_{b}\Gamma,A_{\iota}$.

\item[$(cut)$]
There exist $a_{0}<a$ and $C$
such that $\mbox{{\rm rk}}(C)<b$ and
$(\mathcal{H}_{\gam},\Tht)\vdash^{a_{0}}_{b}\Gamma,\lnot C$
and
$(\mathcal{H}_{\gam},\Tht)\vdash^{a_{0}}_{b}C,\Gamma$.

\item[$(\Ome\in M_{2})$]
There exist ordinals $a_{\ell}$, $a_{r}(\alp)$ and a sentence $C\in\Pi_{2}(\Ome)$
such that $\sup\{a_{\ell}+1,a_{r}(\alp)+1: \alp<\Ome\}\leq a$, $b\geq\Ome$,
$ (\calh_{\gam},\Tht)\vdash^{a_{\ell}}_{b}\Gam,C$
and
$(\calh_{\gam},\Tht\cup\{\ome\alp\})\vdash^{a_{r}(\alp)}_{b}\lnot C^{(\alp,\Ome)}, \Gam$
for any $\alp<\Ome$.

\item[$({\rm rfl}(\pi,k,\vec{\xi},\vec{\nu}))$]
There exist a Mahlo ordinal $\mK\geq\pi\in\calh_{\gam}[\Tht]\cap(b+1)$, an integer $2\leq k\leq N$ and
sequences $\vec{\nu}=(\nu_{2},\ldots,\nu_{N-1}),\vec{\xi}=(\xi_{2},\ldots,\xi_{N-1})\in SD$ of 
ordinals $\nu_{i},\xi_{i}\in E$,
ordinals $a_{\ell},a_{r}(\rho), a_{0}$,  
and a finite set $\Del$ of $\Sig_{k}(\pi)$-sentences enjoying the following conditions:
When $\pi=\mK$, 
$k=N$ and 
$\vec{\nu}=\vec{0}$ with $lh(\vec{\nu})=N-1$ hold.
Also let $\vec{\xi}=\vec{0}$ in this case.
When $\pi<\mK$, 
$\xi_{k}\neq 0$ with $k<N$,
$\vec{0}\neq\vec{\xi}$, and
$\fal i(\xi_{i}\leq_{sp} m_{i}(\pi))$.


\benu

\item
When $\pi<\mK$, cf.\,Definitions \ref{df:Lam}.\ref{df:Exp2.10},
\beqn\label{eq:coeffbnd}
\fal i<k(\nu_{i}=\xi_{i}) \spand (\nu_{k},\ldots,\nu_{N-1})<_{sd}\xi_{k}
\spand K(\vec{\nu})\cup K(\vec{\xi})\subset \calh_{\gam}[\Tht]
\eeqn
and
\begin{equation}\label{eq:notation.12rfl}
\forall\mu\in \vec{\nu}\cup\vec{\xi}\cup \vec{m}(\pi)(K(\mu)\subset \mathcal{H}_{\max K(\mu)}[\Theta])
\end{equation}
cf.\,(\ref{eq:notationsystem.12}).




 \item
For each $\del\in\Del$,
 $(\mathcal{H}_{\gam},\Tht)\vdash^{a_{\ell}}_{b}\Gamma, \lnot\del$.

\item
$H(\vec{\nu},\pi,\gam,\Tht)$ denotes the \textit{resolvent class} for 
$\vec{\nu}$, $\pi$, $\gam$ and $\Tht$ defined as follows:
{\small
\beqnarr 
C(\pi,\gam,\Tht) & := & \{\rho<\pi :\calh_{\gam}(\rho)\cap\pi\subset\rho \spand \Tht\cap\pi\subset\rho\}
\label{eq:Hdfk}
\\
\rho\in H(\vec{\nu},\pi,\gam,\Tht) & :\Leftrightarrow & \fal i(\nu_{i}\leq_{sp} m_{i}(\rho)\land 
K(m_{i}(\rho))\subset\mathcal{H}_{\max K(m_{i}(\rho))}(\rho))
\nonumber
\eeqnarr
}
for $\rho\in Reg\cap C(\pi,\gam,\Tht) $.


Then for each $\rho\in H(\vec{\nu},\pi,\gam,\Tht)$,
$
(\mathcal{H}_{\gam},\Tht\cup\{\rho\})\vdash^{a_{r}(\rho)}_{b}\Gamma, 
\Del^{(\rho,\pi)}
$.

\item
\beqn\label{eq:Mhordgap}
\sup\{a_{\ell},a_{r}(\rho): \rho\in H(\vec{\nu},\pi,\gam,\Tht) \}\leq a_{0}\in\calh_{\gam}[\Tht]\cap a
\eeqn

\eenu

\edes
}
\end{definition}

In the inference rule $({\rm rfl}(\pi,k,\vec{\xi},\vec{\nu}))$ for $\pi=\psi_{\sig}^{\vec{\xi}}(c)<\mK$,
we have $\pi\in Mh^{c}_{2}(\vec{\xi})$. 
In particular, $\pi\in \bigcap_{i<k}Mh^{c}_{i}(\xi_{i}) \cap Mh^{c}_{k}(\xi_{k})$.
Also we are assuming
$(\nu_{k},\ldots,\nu_{N-1})<_{sd}\xi_{k}$, a fortiori $(\nu_{k},\ldots,\nu_{N-1})<\xi_{k}$.
Since $\pi\in  \bigcap_{i<k}Mh^{c}_{i}(\nu_{i})$ is a $\Pi_{k}$-sentence holding on $L_{\pi}$,
we obtain $\pi\in M_{k}(Mh^{c}_{2}(\vec{\nu}))$.
Thus the reflection rule $({\rm rfl}(\pi,k,\vec{\nu}))$ says that
$\pi$ is $\Pi_{k}$-reflecting on the class
$H(\vec{\nu},\pi,\gam,\gam_{0},\Tht)$
for the club subset 
$C(\pi,\gam,\Tht)$ 
of $\pi$, 
cf.\,Proposition \ref{th:KM1}.
On the other side 
we see $\rho\in Mh_{2}^{a}(\vec{\nu})$ from Proposition \ref{prp:Mhelementary}
if $\fal i(\nu_{i}\leq m_{i}(\rho))$ for $\rho\in Mh_{2}^{a}(\vec{m}(\rho))$.

We will state some lemmas for the operator controlled derivations.
These can be shown as in \cite{Buchholz}.
In what follows by an operator $\mathcal{H}$ we mean an $\mathcal{H}_{\gamma}$ for an ordinal $\gamma$.

\blem\label{lem:weakening}
Let $(\calh_{\gam},\Tht)\vdash^{a}_{b}\Gam$.
\benu
\item\label{lem:weakening1}
$(\calh_{\gam'},\Tht\cup\Tht_{0})\vdash^{a'}_{b'}\Gam,\Del$ for any $\gam'\geq\gam$,
any $\Tht_{0}$, and any $a'\geq a$, $b'\geq b$ such that
$\sfk(\Del)\cup\{a'\}\subset\calh_{\gam'}[\Tht\cup\Tht_{0}]$.

\item\label{lem:weakening2}
Assume $\Tht_{1}\cup\{c\}=\Tht$,
$c\in\calh_{\gam}[\Tht_{1}]$.
Then $(\calh_{\gam},\Tht_{1})\vdash^{a}_{b}\Gam$.
\eenu
\elem

\begin{lemma}\label{lem:tautology}{\rm (Tautology)}
$(\mathcal{H},{\sf k}(\Gamma\cup\{A\}))\vdash^{2\footnotesize{\mbox{{\rm rk}}}(A)}_{0}\Gamma,\lnot A, A$.
\end{lemma}

\blem\label{lem:inversionreg}{\rm (Inversion)}
Let  $A\simeq \bigwedge(A_{\iota})_{\iota\in J}$, and
$(\calh,\Tht)\vdash^{a}_{b}\Gam$ with $A\in\Gam$.
Then for any $\iota\in J$, 
$(\calh,\Tht\cup{\sf k}(\iota))\vdash^{a}_{b}\Gam,A_{\iota}$ holds.
\elem

\begin{lemma}\label{lem:boundednessreg}{\rm (Boundedness)}
Suppose $(\mathcal{H},\Tht)\vdash^{a}_{c}\Gam, C$ for a $C\in\Sig_{1}(\lam)$,
{\rm and} $a\leq b\in\mathcal{H}\cap\lambda$. Then
$(\mathcal{H},\Tht)\vdash^{a}_{c}\Gam,C^{(b,\lam)}$.

\end{lemma}

\blem\label{lem:persistency}{\rm (Persistency)}
Suppose $(\mathcal{H},\Tht)\vdash^{a}_{c}\Gam, C^{(b,\lam)}$ for a $C\in\Sig_{1}(\lam)$
and a $b<\lambda\in\calh[\Tht]$. Then
$(\mathcal{H},\Tht)\vdash^{a}_{c}\Gam,C$.
\elem

\begin{lemma}\label{lem:predcereg}{\rm (Predicative Cut-elimination)}
Suppose $(\mathcal{H},\Tht)\vdash^{b}_{c+\omega^{a}}\Gamma$, $a\in\mathcal{H}[\Tht]$
and $]c,c+\ome^{a}]\cap Reg=\emptyset$.
Then $(\mathcal{H},\Tht)\vdash^{\varphi ab}_{c}\Gamma$.

\end{lemma}

\begin{lemma}\label{th:embedreg}{\rm (Embedding of Axioms)}\\
 For each axiom $A$ in $\mbox{{\sf KP}}\Pi_{N}$, there is an $m<\omega$ such that
 for any operator $\mathcal{H}=\calh_{\gam}$,  
 $(\mathcal{H},\emptyset)\vdash^{\mK\cdot 2}_{\mK+m}  A$
holds.
\end{lemma}
\bprf
The axiom $\lnot A,\exi z\, A^{(z)}$ for $\Pi_{N}$-reflection follows from 
$A,\lnot A$ and 
\\
$\exi z\, A^{(z)},\lnot A^{(\rho)}$ for regular ordinals $\rho<\mK$
by an inference $({\rm rfl}(\mK,N,\vec{0},\vec{0}))$.
\eprf

\blem\label{th:embedregthm}{\rm (Embedding)}
If $\mbox{{\sf KP}}\Pi_{N}\vdash \Gam$ for sets $\Gam$ of sentences, 
there are $m,k<\ome$ such that for any operator $\calh=\calh_{\gam}$, 
 $(\calh,\emptyset)\vdash_{\mK+m}^{\mK\cdot 2+k}\Gam$ 
holds
\elem

\section{Lowering and eliminating higher Mahlo operations}\label{subsec:elimpi11}
In the section inferences $({\rm rfl}(\mK,N,\vec{0},\vec{0}))$
for $\Pi_{N}$-reflecting ordinals $\mK$ are eliminated from operator controlled derivations of $\Sig_{1}$-sentences 
$\varphi^{L_{\Ome}}$ over $\Ome$.


$\alp\#\bet$ denotes the natural (commutative) sum of ordinal terms $\alp,\bet$.

\blem\label{lem:KppiNlower}
For a Mahlo term $\pi\in OT$, $\vec{\xi}\in SD$ 
denotes a sequence with $lh(\vec{\xi})=N-2$, and $2\leq k\leq N-1$ an integer for which the following hold:
When $\pi=\mK$, let $\vec{\xi}=\vec{0}$ and $k=N-1$.
Otherwise 
$\vec{\xi}=(\xi_{2},\ldots,\xi_{k+1})*\vec{0}$ with
$\xi_{k+1}\neq 0$ such that $\fal i\leq k+1(\xi_{i}\leq_{sp}m_{i}(\pi))$.

For ordinal terms $\gam,a\in OT$ 
let us define a sequence
$\vec{\zeta}(a)  :=  (\zeta_{2}(a),\ldots,\zeta_{k}(a))*\vec{0}$ with $lh(\vec{\zeta}(a))=N-2$ as follows.
$\vec{\zeta}(a)=\vec{0}*(\gam+a)$
when $\pi=\mK$.
Otherwise $\zeta_{k}(a)=\xi_{k}+\Lam^{\xi_{k+1}}(\gam+a)$
and $\zeta_{i}(a)=\xi_{i}$ for $i<k$.

Let $\kap\in H(\vec{\zeta}(a),\pi,\gam,\Tht)$ for a finite set $\Tht\subset OT$.

Now suppose 
$
(\mathcal{H}_{\gam},\Tht)\vdash^{a}_{\pi}\Gamma
$
where $\{\gam,\pi\}\cup K(\vec{\xi})\subset\calh_{\gam}[\Tht]$, $\Tht\subset\pi$,
$\fal i(K(\xi_{i})\subset\mathcal{H}_{\max K(\xi_{i})}[\Tht])$,
and $\Gamma\subset\Pi_{k+1}(\pi)$.

Let
$\gam(a,b)=\gam\#a\#b$, $\bet(a,b)=\psi_{\pi}(\gam(a,b))$, and $c>\gam(a,\kap)$.
Then
the following holds:
\beqn\label{eq:KppiNlowerconcl}
(\mathcal{H}_{c},\Tht\cup\{\kap\})\vdash^{\bet(a,\kap)}_{\kap}
\Gamma^{(\kap,\pi)}
\eeqn

\elem
{\bf Proof} by induction on $a$. 
Let $\kap\in H(\vec{\zeta}(a),\pi,\gam,\Tht)$.
We see $\vec{\zeta}(a)\in SD$, and from (\ref{eq:controlder}) and $\Tht\subset\kap$
that
\beqn\label{eq:CollapsingthmKR100}
\sfk(\Gam)\cap\pi\subset\calh_{\gam}(\kap)\cap\pi\subset \kap
\eeqn

For any $a\in\calh_{\gam}[\Tht]$, we obtain $\{\gam,\pi,a,\kap\}\subset\calh_{\gam}(\pi)$ by 
$\Tht\cup\{\kap\}\subset\pi$.
Hence for $\gam(a,\kap)=\gam\#a\#\kap$, $\{\gam(a,\kap),\pi\}\subset\calh_{\gam}(\pi)$, and 
$\{\gam(a,\kap),\pi\}\subset\calh_{\gam(a,\kap)}(\bet(a,\kap))$
by the definition (\ref{eq:Psivec}).
Therefore $\kap\in\calh_{\gam(a,\kap)}(\bet(a,\kap))\cap\pi\subset\bet(a,\kap)$ by Proposition \ref{prp:Cantorclosed},
and $\Tht\subset\bet(a,\kap)<\pi$.
Thus we obtain
\[
\{a_{0},a_{1}\}\subset\calh_{\gam}[\Tht\cup\Tht_{0}]\spand a_{0}<a_{1} \spand \Tht_{0}\subset\kap
\Rarw \bet(a_{0},\kap)<\bet(a_{1},\kap).
\]
{\bf Case 1}. 
First consider the case when the last inference is a $({\rm rfl}(\pi,k+1,\vec{\xi},\vec{\nu}))$.

We have $a_{\ell}\in\mathcal{H}_{\gam}[\Tht]\cap a$, 
$a_{r}(\rho)\in\mathcal{H}_{\gam}[\Tht\cup\{\rho\}]\cap a$, and 
a finite set $\Del$ of $\Sig_{k+1}(\pi)$-sentences.
We have for each $\del\in\Del$
\beqn\label{eq:KppiNlowerCase1a}
(\mathcal{H}_{\gam},\Theta)\vdash^{a_{\ell}}_{\pi}\Gamma,\lnot\del
\eeqn
and for each $\rho\in H(\vec{\nu},\pi,\gam,\Tht)$
\beqn\label{eq:KppiNlowerCase1b}
 (\mathcal{H}_{\gam},\Theta\cup\{\rho\})\vdash^{a_{r}(\rho)}_{\pi}
 \Gamma, \Del^{(\rho,\pi)}
\eeqn
 
When $\pi<\mK$,
 $\vec{\nu}=(\nu_{2},\ldots,\nu_{N-1})\in SD$ is a sequence 
 such that
 $\fal i<k+1(\nu_{i}=\xi_{i})$, $(\nu_{k+1},\ldots,\nu_{N -1})<_{sd}\xi_{k+1}$,
$K(\vec{\nu})\cup K(\vec{\xi})\subset\calh_{\gam}[\Tht]$, 
and
$\fal i(K(\nu_{i})\subset \mathcal{H}_{\max K(\nu_{i})}[\Theta])$,
cf.\,(\ref{eq:coeffbnd}) and (\ref{eq:notation.12rfl}).

Let $\Gam_{0}=\Gam\cap\Sig_{k}(\pi)$ and $\{\fal x\in L_{\pi}\,\tht_{i}(x): i=1,\ldots,n\}\,(n\geq 0)=\Gam\setm\Gam_{0}$
for $\Sig_{k}(\pi)$-formulas $\tht_{i}(x)$.
Let us fix $\vec{d}=\{d_{1},\ldots,d_{n}\}\subset Tm(\kap)$ arbitrarily.
Put ${\sf k}(\vec{d})=\bigcup\{{\sf k}(d_{i}) : i=1,\ldots,n\}$ and $\Gam(\vec{d})=\Gam_{0}\cup\{\tht_{i}(d_{i}): i=1,\ldots,n\}$.

By Inversion lemma \ref{lem:inversionreg} from (\ref{eq:KppiNlowerCase1a}) we obtain for  each $\del\in\Del$
\beqn\label{eq:L4.10.N}
(\mathcal{H}_{\gam},\Theta\cup{\sf k}(\vec{d}))\vdash^{a_{\ell}}_{\pi}\Gamma(\vec{d}),\lnot\del
\eeqn

Let $\rho\in C(\kap,c,\Tht\cup\{\kap\}\cup{\sf k}(\vec{d}))$.
We see $\rho<\kap$, and ${\sf k}(\vec{d})<\rho$
from ${\sf k}(\vec{d})<\kap$.
By $\Tht\cap\pi\subset\calh_{\gam}(\kap)\cap\pi\subset\kap$ and $\gam\leq c$ 
we obtain $C(\kap,c,\Tht\cup\{\kap\}\cup{\sf k}(\vec{d}))\subset C(\pi,\gam,\Tht)$. 
Namely, cf.\,(\ref{eq:Hdfk})
\beqn\label{eq:L4.10.0}
\rho\in H(\vec{\nu},\kap,c,\Tht\cup\{\kap\}\cup{\sf k}(\vec{d}))
\Rarw \rho\in H(\vec{\nu},\pi,\gam,\Tht)
\eeqn
For each $\rho\in H(\vec{\nu},\kap,c,\Tht\cup\{\kap\}\cup{\sf k}(\vec{d}))$,
IH with (\ref{eq:KppiNlowerCase1b}) and (\ref{eq:L4.10.0})
yields for $c>\gam(a_{r}(\rho),\kap)$
and $\kap\in H(\vec{\zeta}(a_{r}(\rho)),\pi,\gam,\Tht\cup\{\rho\})$
\beqn\label{eq:L4.10case1r}
(\calh_{c},\Tht\cup\{\rho,\kap\})
\vdash^{\bet(a_{r}(\rho),\kap)}_{\kap}
\Gamma^{(\kap,\pi)},\Del^{(\rho,\pi)}
\eeqn
Let
$\rho\in M_{\ell} 
:=
\{\rho\in Reg: \fal i(\zeta_{i}(a_{\ell})\leq_{sp} m_{i}(\rho))\} \cap H(\vec{\nu},\kap,c,\Tht\cup\{\kap\}\cup{\sf k}(\vec{d}))
$.
Then 
$M_{\ell}\subset H(\vec{\zeta}(a_{\ell}),\pi,\gam,\Tht\cup{\sf k}(\vec{d}))$ and
$\Tht\cup{\sf k}(\vec{d})\subset\rho$.
For each $\del\in\Del$, 
IH with (\ref{eq:L4.10.N})  yields for $c>\gam(a_{\ell},\rho)$
\beqn\label{eq:L4.10case1l}
(\calh_{c},\Tht\cup{\sf k}(\vec{d})\cup\{\rho\})\vdash^{\bet(a_{\ell},\rho)}_{\rho}
\Gamma(\vec{d})^{(\rho,\pi)},\lnot\del^{(\rho,\pi)}
\eeqn
From (\ref{eq:L4.10case1r}) and (\ref{eq:L4.10case1l}) by several $(cut)$'s of $\del^{(\rho,\pi)}$ 
with $\mbox{rk}(\del^{(\rho,\pi)})<\kap$ we obtain for $a(\rho)=\max\{a_{\ell},a_{r}(\rho)\}$
and some $p<\ome$
\beqn \label{eq:L4.10case1.1a}
\{(\calh_{c},\Tht\cup{\sf k}(\vec{d})\cup\{\kap,\rho\})
\vdash^{\bet(a(\rho),\kap)+p}_{\kap}
\Gamma(\vec{d})^{(\rho,\pi)} , \Gamma^{(\kap,\pi)}
: 
\rho\in M_{\ell}\}
\eeqn
On the other hand
we have by Tautology lemma \ref{lem:tautology}
for each $ \tht(\vec{d})^{(\kap,\pi)}\in\Gam(\vec{d})^{(\kap,\pi)}$
\beqn\label{eq:L4.10case1.1b}
(\calh_{\gam},\Tht\cup{\sf k}(\vec{d})\cup\{\kap\})\vdash^{2\footnotesize{\mbox{{\rm rk}}}(\tht(\vec{d})^{(\kap,\pi)})}_{0} 
\Gam(\vec{d})^{(\kap,\pi)}, \lnot\tht(\vec{d})^{(\kap,\pi)}
\eeqn
where 
$2\mbox{{\rm rk}}(\tht(\vec{d})^{(\kap,\pi)}) \leq\kap+p$ for some $p<\ome$.

Moreover we have
$\sup\{2\mbox{{\rm rk}}(\tht(\vec{d})^{(\kap,\pi)}), \bet(a(\rho),\kap)+p: \rho\in M_{\ell}\}
\leq \bet(a_{0},\kap)+p\in\calh_{\gam}[\Tht\cup\{\kap\}]$,
where $\sup\{ a_{\ell}, a_{r}(\rho): \rho\in H(\vec{\nu},\pi,\gam,\Tht)\}\leq a_{0}<a$  
 by (\ref{eq:Mhordgap}).
 
Now let $\vec{\mu}=(\mu_{2},\ldots,\mu_{N-1})=\max\{\vec{\zeta}(a_{\ell}),\vec{\nu}\}$ with
$\mu_{i}=\max\{\zeta_{i}(a_{\ell}),\nu_{i}\}$.
Since $\nu_{i}=\xi_{i}\leq_{pt}\zeta_{i}(a_{\ell})$ for $i<k+1$, we obtain
$\mu_{i}=\left\{
\begin{array}{ll}
\zeta_{i}(a_{\ell}) & i\leq k
\\
\nu_{i} & i>k
\end{array}
\right.
$.
We see that
$M_{\ell}=H(\vec{\mu},\kap,c,\Tht\cup\{\kap\}\cup{\sf k}(\vec{d}))$.
Moreover we have $\fal i<k(\mu_{i}=\xi_{i}=\zeta_{i}(a))$ and 
$(\mu_{k},\ldots,\mu_{N-1})=(\zeta_{k}(a_{\ell}))*(\nu_{k+1},\ldots,\nu_{N -1})<_{sd}\zeta_{k}(a)$.
Also
$\forall i(K(\zeta_{i}(a))\subset\mathcal{H}_{\max K(\zeta_{i}(a))}[\Theta])$ and $\forall i(K(\mu_{i})\subset\mathcal{H}_{\max K(\mu_{i})}[\Theta])$.
For $\lnot\Gam(\vec{d})^{(\kap,\pi)}\subset\Pi_{k}(\kap)$, 
by an inference rule 
$({\rm rfl}(\kap,k,\vec{\zeta}(a),\vec{\mu}))$ with 
its resolvent class 
$M_{\ell}$,
we conclude
from (\ref{eq:L4.10case1.1b})
and (\ref{eq:L4.10case1.1a}) that
$
(\calh_{c},\Tht\cup\{\kap\}\cup{\sf k}(\vec{d}))
\vdash^{\bet(a_{0},\kap)+p+1}_{\kap}
\Gam(\vec{d})^{(\kap,\pi)}, \Gamma^{(\kap,\pi)}
$.
Since $\vec{d}\subset Tm(\kap)$ is arbitrary,
 several $(\bigwedge)$'s yield (\ref{eq:KppiNlowerconcl}).
\\

\noindent
{\bf Case 2}.
Second consider the case when the last inference is a $({\rm rfl}(\pi,j,\vec{\xi},\vec{\nu}))$
for a $j<k+1$.
We have 
$(\mathcal{H}_{\gam},\Theta)\vdash^{a_{\ell}}_{\pi}\Gamma,\lnot\del$
for each $\del\in\Del\subset\Sig_{j}(\pi)$ with $a_{\ell}\in\mathcal{H}_{\gam}[\Tht]\cap a$,
and
$(\mathcal{H}_{\gam},\Theta\cup\{\rho\})\vdash^{a_{r}(\rho)}_{\pi} \Gamma,\Del^{(\rho,\pi)}$
for each 
$\rho\in H(\vec{\nu},\pi,\gam,\Tht)$ with $a_{r}(\rho)\in\mathcal{H}_{\gam}[\Tht\cup\{\rho\}]\cap a$.
$\vec{\nu}\in SD$ is a sequence such that 
$\fal i<j(\nu_{i}=\xi_{i})$ and $(\nu_{j},\ldots,\nu_{N-1})<_{sd}\xi_{j}$.

We see that the resolvent class $H(\vec{\nu},\kap,c_{1},\Tht\cup\{\kap\})$ is a subclass of $H(\vec{\nu},\pi,\gam,\Tht)$.
By IH we have 
$
(\mathcal{H}_{c},\Theta\cup\{\kap\})
\vdash^{\bet(a_{\ell},\kap)}_{\kap}
\Gam^{(\kap,\pi)},\lnot\del^{(\kap,\pi)}
$
for each $\del\in\Del$,
and 
$
(\mathcal{H}_{c},\Theta\cup\{\kap,\rho\} )
\vdash^{\bet(a_{r}(\rho),\kap)}_{\kap}
 \Gamma^{(\kap,\pi)},\Del^{(\rho,\pi)}
$
for each $\rho\in H(\vec{\nu},\kap,c,\Tht\cup\{\kap\})$ 
with $\Del^{(\rho,\pi)}=(\Del^{(\kap,\pi)})^{(\rho,\kap)}$.
We claim that $\fal i\leq j(\xi_{j}\leq_{sp}m_{i}(\kap))$.
Consider the case when $i=j=k$.
Then we have $\xi_{k}\leq_{sp}m_{k}(\pi)$ and $\zeta_{k}(a)\leq_{sp}m_{k}(\kappa)$
with $\xi_{k}<_{pt}\zeta_{k}(a)$.
We obtain $\xi_{k}\leq_{sp}m_{k}(\kappa)$.
Hence by an inference rule $({\rm rfl}(\kap,j,\vec{\xi}(j),\vec{\nu}))$ 
for the sequence $\vec{\xi}(j)=(\xi_{2},\ldots,\xi_{j})*\vec{0}\in SD$, cf.\,Proposition \ref{prp:SD}.\ref{prp:SD.-1},
 we obtain (\ref{eq:KppiNlowerconcl}).
\\
 
 \noindent
{\bf Case 3}.
Third consider the case when the last inference is a $({\rm rfl}(\sig,j,\vec{\mu},\vec{\nu}))$
for a $\sig<\pi$.
We have 
$(\mathcal{H}_{\gam},\Tht)\vdash^{a_{\ell}}_{\pi}\Gamma,\lnot\del$
for each $\del\in\Del\subset\Sig_{j}(\sig)$,
and 
$(\mathcal{H}_{\gam},\Tht\cup\{\rho\} )\vdash^{a_{r}(\rho)}_{\pi}
\Gamma,\Del^{(\rho,\sig)}$
for each $\rho\in H(\vec{\nu},\sig,\gam,\Tht)$.
We obtain $\sig<\kap$
by (\ref{eq:CollapsingthmKR100}) for $\sig\in\calh_{\gam}[\Tht]$.
Hence
$\Del\subset\Sig^{1}_{0}(\sig)\subset\Sig_{0}(\kap)$ and 
$\del^{(\kap,\pi)}\equiv \del$ for any $\del\in\Del$.
Let $H(\vec{\nu},\sig,c,\Tht\cup\{\kap\})$ be the resolvent class for $\sig$, $\vec{\nu}$, $c$ and $\Tht\cup\{\kap\}$.
Then $H(\vec{\nu},\sig,c,\Tht\cup\{\kap\})\subset H(\vec{\nu},\sig,\gam,\Tht)$.

From IH
we have
$(\mathcal{H}_{c},\Tht\cup\{\kap\})\vdash^{\bet(a_{\ell},\kap)}_{\kap}\Gamma^{(\kap,\pi)},\lnot\del$ 
for each $\del\in\Del$, and
$(\mathcal{H}_{c},\Tht\cup\{\kap,\rho\})
\vdash^{\bet(a_{r}(\rho),\kap)}_{\kap}\Gamma^{(\kap,\pi)},\Del^{(\rho,\sig)}$ 
for each $ \rho\in H(\vec{\nu},\sig,c,\Tht\cup\{\kap\})$.
 We obtain 
(\ref{eq:KppiNlowerconcl})
 by an inference rule $({\rm rfl}(\sig,j,\vec{\mu},\vec{\nu}))$
with the resolvent class $H(\vec{\nu},\sig,c,\Tht\cup\{\kap\})$.
\\

\noindent
{\bf Case 4}. 
Fourth consider the case when the last inference $(\bigwedge)$ introduces a $\Pi_{k+1}(\pi)$-sentence 
$(\fal x\in L_{\pi}\,\tht(x))\in\Gam$.
We have 
$(\mathcal{H}_{\gam},\Theta\cup{\sf k}(d))\vdash^{a(d)}_{\pi}\Gamma,\tht(d)$
for each $d\in Tm(\pi)$.
For each $d\in Tm(\kap)$, 
IH with ${\sf k}(d)<\kap$ yields
$(\mathcal{H}_{c},\Theta\cup\{\kap\}\cup{\sf k}(d))\vdash^{\bet(a(d),\kap)}_{\kap}
\Gamma^{(\kap,\pi)},
\tht(d)^{(\kap,\pi)}$.
$(\bigwedge)$ yields (\ref{eq:KppiNlowerconcl}) 
for $\fal x\in L_{\kap}\,\tht(x)^{(\kap,\pi)}\equiv(\fal x\in L_{\pi}\,\tht(x))^{(\kap,\pi)}\in\Gamma^{(\kap,\pi)}$.
\\

\noindent
{\bf Case 5}. 
Fifth consider the case when the last inference $(\bigwedge)$ introduces a $\Sig_{0}(\pi)$-sentence 
$(\fal x\in c\,\tht(x))\in\Gam$ for a $c\in Tm(\pi)$.
We have 
$(\mathcal{H}_{\gam},\Theta\cup{\sf k}(d))\vdash^{a(d)}_{\pi}\Gamma,\tht(d)$
for each $d\in Tm(|c|)$.
Then we have $|d|<|c|<\kap$ by (\ref{eq:CollapsingthmKR100}).
IH yields $(\mathcal{H}_{c},\Theta\cup\{\kap\}\cup{\sf k}(d)\vdash^{\bet(a(d),\kap)}_{\kap}
\Gamma^{(\kap,\pi)},\tht(d)$, and we obtain (\ref{eq:KppiNlowerconcl}) by an inference $(\bigwedge)$.
\\

\noindent
{\bf Case 6}.
Sixth consider the case when the last inference $(\bigvee)$ introduces a $\Sig_{k}(\pi)$-sentence 
$(\exi x\in L_{\pi}\,\tht(x))\in\Gam$.
We have 
$(\calh_{\gam},\Tht)\vdash^{a_{0}}_{\pi}\Gam,\tht(d)$
for a $d\in Tm(\pi)$.
Without loss of generality we can assume that ${\sf k}(d)\subset {\sf k}(\tht(d))$.
Then we see that $|d|<\kap$ from (\ref{eq:CollapsingthmKR100}), and $d\in Tm(\kap)$.
Also $|d|<\kap<\bet(a,\kap)$ for (\ref{eq:bigveebnd}).
IH yields with $(\exi x\in L_{\pi}\,\tht(x))^{(\kap,\pi)}\equiv(\exi x\in L_{\kap}\, \tht(x)^{(\kap,\pi)})\in\Gam^{(\kap,\pi)}$,
$(\calh_{c},\Tht\cup\{\kap\})\vdash^{\bet(a_{0},\kap)}_{\kap}\Gam^{(\kap,\pi)},\tht(d)^{(\kap,\pi)}$,
and we obtain (\ref{eq:KppiNlowerconcl}) by an inference $(\bigvee)$.
\\

\noindent
{\bf Case 7}.
Seventh consider the case when the last inference is a $(cut)$.
We have $(\calh_{\gam},\Tht)\vdash^{a_{0}}_{\pi}\Gam,\lnot C$ and
$(\calh_{\gam},\Tht)\vdash^{a_{0}}_{\pi}C,\Gam$ for $a_{0}<a$ with $\mbox{rk}(C)<\pi$.
Then $C\in\Sig_{0}(\pi)$ by Proposition \ref{lem:rank}.\ref{lem:rank3}.
On the other side ${\sf k}(C)\subset\pi$ holds
by Proposition \ref{lem:rank}.\ref{lem:rank1}.
Then ${\sf k}(C)\subset\kap$
by (\ref{eq:CollapsingthmKR100}).
Hence $C^{(\kap,\pi)}\equiv C$ and
 $\mbox{rk}(C^{(\kap,\pi)})<\kap$ again by Proposition \ref{lem:rank}.\ref{lem:rank1}.
IH yields
$ (\calh_{c},\Tht\cup\{\kap\})\vdash^{\bet(a_{0},\kap)}_{\kap}\Gam^{(\kap,\pi)},\lnot C^{(\kap,\pi)}$
and
 $(\calh_{c},\Tht\cup\{\kap\})\vdash^{\bet(a_{0},\kap)}_{\kap}C^{(\kap,\pi)},\Gam^{(\kap,\pi)}$.
 Hence by a $(cut)$ we obtain (\ref{eq:KppiNlowerconcl}).
\\

\noindent
{\bf Case 8}.
Eighth consider the case when the last inference is an $(\Ome\in M_{2})$.
We have $(\calh_{\gam},\Tht)\vdash^{a_{\ell}}_{\pi}\Gam,C$ and
$(\calh_{\gam},\Tht\cup\{\ome\alp\})\vdash^{a_{r}(\alp)}_{\pi}\lnot C^{(\alp,\Ome)}, \Gam$
for each $\alp<\Ome$
with $\sup\{a_{\ell}+1,a_{r}(\alp)+1: \alp<\Ome\}\leq a$ and $C\in\Pi_{2}(\Ome)$.

We obtain $\ome\alp<\kap$ for $\alp<\Ome$.
IH with $C^{(\kap,\pi)}\equiv C$ yields  for each $\alp<\Ome$,
$(\calh_{c},\Tht\cup\{\kap,\ome\alp\})\vdash^{\bet(a_{r}(\alp),\kap)}_{\kap}\lnot C^{(\alp,\Ome)}, \Gam^{(\kap,\pi)}$,
and
$ (\calh_{c},\Tht\cup\{\kap\})\vdash^{\bet(a_{\ell},\kap)}_{\kap}\Gam^{(\kap,\pi)},C$.
 An $(\Ome\in M_{2})$ yields
(\ref{eq:KppiNlowerconcl})

All other cases are seen easily from IH.
\eprf

\blem\label{lem:collapsKN}
Let  $\lam\leq\pi$ be a regular ordinal term such that
$\fal i(K(m_{i}(\pi))\subset\mathcal{H}_{\max K(m_{i}(\pi))}[\Tht])$,
and $\Gamma\subset\Sig_{1}(\lam)$.

Suppose for an ordinal term $a\in OT$ 
\[
(\mathcal{H}_{\gam},\Tht)\vdash^{a}_{\pi}\Gamma
\]
where $\{\gam,\lam,\pi\}\subset\calh_{\gam}[\Tht]$.



Assume
\beqn\label{eq:asscollasKN}
\fal\rho\in[\lam,\pi] \fal d[\Tht\subset\psi_{\rho}(\gam\#d)]
\eeqn

Let
$\hat{a}=\gam\# \ome^{\pi+a+1}$ and $\bet=\psi_{\lam}(\hat{a})$.
Then 
the following holds
\beqn\label{eq::collapsKNconcl}
(\mathcal{H}_{\hat{a}+1},\Tht)
\vdash^{\bet}_{\bet}
\Gamma
\eeqn

\elem
{\bf Proof} by main induction on $\pi$ with subsidiary induction on $a$.
We can assume $a>0$.



We see that $\Tht\subset\bet=\psi_{\lam}(\hat{a})$ from (\ref{eq:asscollasKN}).
Hence 
\[
a_{0}\in\calh_{\gam}[\Tht]\cap a  \Rarw \psi_{\lam}(\widehat{a_{0}})<\psi_{\lam}(\hat{a})
\]

Let $\vec{\xi}\in SD$ be a sequence of ordinals and $k$ a number for which the following hold:
If $\pi=\mK$, then let $\vec{\xi}=\vec{0}$ with $lh(\vec{\xi})=N-1$ and $k=N-1$.
Let $\pi<\mK$. If $\vec{m}(\pi)\neq\vec{0}$, then $K(\vec{\xi})\subset\calh_{\gam}[\Tht]$,
$\vec{\xi}\leq \vec{m}(\pi)$ and $k=\max\{k\leq N-2: \xi_{k+1}>0\}$.
Otherwise let $\vec{\xi}=\vec{0}$ and $k=1$. 
By the assumption (\ref{eq:asscollasKN}), 
and (\ref{eq:controlder}) we obtain
\beqn\label{eq:collapsKNclose}
\fal\rho\in[\lam,\pi] \fal b\in K(\vec{\xi})\fal d
[{\sf k}(\Gam)\cup\{\gam,\lam,a,\pi,b\} \subset
\calh_{\gam}(\psi_{\rho}(\gam\#d))]
\eeqn
{\bf Case 1}. 
First consider the case when $k\geq2$. 

Let $\vec{\xi}=\vec{m}(\pi)$,
and $\vec{\zeta}(a):= (\zeta_{2}(a),\ldots,\zeta_{k}(a))*\vec{0}$ be the 
sequence defined as in Lemma \ref{lem:KppiNlower} from $\gam,a$:
 $\vec{\zeta}(a)=\vec{0}*(\gam+a)$ when $\pi=\mK$, otherwise
$\zeta_{k}(a)=\xi_{k}+\Lam^{\xi_{k+1}}(\gam+a)$
and $\zeta_{i}(a)=\xi_{i}$ for $i<k$.
Also let $\gam(a,b)=\gam\#a\#b$ and $\bet(a,b)=\psi_{\pi}\gam(a,b)$.



Let $\kap:=\psi_{\pi}^{\vec{\zeta}(a)}(\gam(a,0))$.
By the assumption (\ref{eq:asscollasKN}) 
we have 
$\Tht\subset\psi_{\pi}(\gam\#a)$.
On the other hand we have $\psi_{\pi}(\gam\#a)=\psi_{\pi}(\gam(a,0))\leq\kap$, and $\Tht\subset\kap$.
$\pi\in\calh_{\gam}[\Tht]$ with $\Tht\subset\pi$ yields
$K(\vec{\xi})=K(\vec{m}(\pi))\subset\calh_{\gam}[\Tht]\subset\calh_{\gam(a,0)}(\kap)$.
Hence
$K(\vec{\xi})\cup\{\pi,\gam(a,0)\}\subset\calh_{\gam(a,0)}(\kap)$, and
$\kap\in OT$ by $\gam(a,0)=\gam\# a>0$ and 
Definition \ref{df:notationsystem}.\ref{df:notationsystem.11} such that
$\kap<\pi$ and $\calh_{\gam}(\kap)\cap\pi\subset\kap$.
Moreover we have
$\fal i(K(\zeta_{i}(a))\subset\mathcal{H}_{\max K(\zeta_{i}(a))}[\Tht])$
by $\fal i(K(m_{i}(\pi))\subset\mathcal{H}_{\max K(m_{i}(\pi))}[\Tht])$ and $\{\gam,a\}\subset\mathcal{H}_{\gam}[\Theta]$ with $\Theta\subset\kappa$.
In other words, $\kap\in H(\vec{\zeta}(a),\pi,\gam,\Tht)$.





By Lemma \ref{lem:KppiNlower} we obtain 
$(\mathcal{H}_{\gam(a,\kap)+1},\Tht\cup\{\kap\})\vdash^{\bet(a,\kap)}_{\kap}
\Gamma^{(\kap,\pi)}$,
and Lemma \ref{lem:weakening}.\ref{lem:weakening2} with $\kap\in\mathcal{H}_{\gam(a,0)+1}[\Tht]$
\beqn\label{eq:collapsKN1}
(\mathcal{H}_{\gam(a,\kap)+1},\Tht)\vdash^{\bet(a,\kap)}_{\kap}
\Gamma^{(\kap,\pi)}
\eeqn

If $\lam=\pi$, then $\Gam^{(\kap,\pi)}\subset\Sig_{1}(\kap)\subset\Sig_{0}(\lam)$.
We have $\Tht\subset\psi_{\pi}(\hat{a})=\bet$, and $\kap\in\calh_{\hat{a}}(\bet)$.
Hence $\{\gam,\pi,a,\kap\}\subset\calh_{\hat{a}}(\bet)$, and 
$\gam(a,\kap)=\gam\#a\#\kap<\gam\#\ome^{\pi+a+1}=\hat{a}$.
Therefore $\kap<\bet(a,\kap)\leq\psi_{\pi}(\hat{a})=\bet$.
We obtain (\ref{eq::collapsKNconcl}) by Persistency lemma \ref{lem:persistency}.

Next consider the case when $\lam<\pi$.
Then $\lam<\kap$ and $\Gam^{(\kap,\pi)}=\Gam$.
We have for (\ref{eq:asscollasKN}),
$\fal d\fal\rho\in[\lam,\kap](\Tht\subset\psi_{\rho}(\gam(a,\kap)+1\#d))$.
By MIH on (\ref{eq:collapsKN1})
we obtain
$
(\mathcal{H}_{b_{0}+1},\Tht)\vdash^{\bet_{0}}_{\bet_{0}}
\Gamma
$
for $\bet_{0}=\psi_{\lam}(b_{0})$ with $b_{0}=(\gam(a,\kap)+1)\#\ome^{\kap+\bet(a,\kap)+1}$.
We have $b_{0}=\gam\#a\#\kap\#1\#\ome^{\bet(a,\kap)+1}<\gam\#\ome^{\pi+a+1}=\hat{a}$
by $\bet(a,\kap)<\pi$.
This yields
$\psi_{\lam}(b_{0})=\bet_{0}<\bet=\psi_{\lam}(\hat{a})$ by $\Tht\subset\bet$ and
$\{\gam,\kap,\pi,a\}\subset\calh_{\hat{a}}(\bet)$.
Hence (\ref{eq::collapsKNconcl}) follows.

In what follows suppose $k=1$.
\\

\noindent
{\bf Case 2}. 
Consider the case when the last inference rule is a $({\rm rfl}(\pi,2,\vec{\xi},\vec{\nu}))$.

We have an ordinal term
$a_{\ell}\in\mathcal{H}_{\gam}[\Tht]\cap a$, and a finite set $\Del$ of $\Sig_{2}(\pi)$-sentences
for which
$(\mathcal{H}_{\gam},\Theta)\vdash^{a_{\ell}}_{\pi}\Gamma,\lnot\del$
holds for each $\del\in\Del$.
On the other hand we have sequences $\vec{\nu},(\xi_{2})*\vec{0}\in SD$ such that
 $\vec{\nu}<_{sd}\xi_{2}$ and $K(\vec{\nu})\cup K(\vec{\xi})\subset \calh_{\gam}[\Tht]$ by (\ref{eq:coeffbnd}), and
an ordinal term
$a_{r}(\rho)\in\mathcal{H}_{\gam}[\Tht\cup\{\rho\}]\cap a$
for which
$(\mathcal{H}_{\gam},\Theta\cup\{\rho\})\vdash^{a_{r}(\rho)}_{\pi}\Gamma, \Del^{(\rho,\pi)}$
holds for each $\rho\in H(\vec{\nu},\pi,\gam,\Tht)$, where 
$\xi_{2}\leq_{sp} m_{2}(\pi)$.



 Let 
  $\rho:=\psi_{\pi}^{\vec{\nu}}(\widehat{a_{\ell}}\#\pi)$
  for
  $\widehat{a_{\ell}}=\gam\#\ome^{\pi+a_{\ell}+1}$.
By the assumption (\ref{eq:asscollasKN}) we have 
$\Tht\subset\psi_{\pi}(\widehat{a_{\ell}})\subset\rho$.
$K(\vec{\nu})\cup\{\pi,\gam,a\}\subset\calh_{\gam}[\Tht]$
yields 
$K(\vec{\nu})\cup\{\pi,\widehat{a_{\ell}}\}\subset\calh_{\widehat{a_{\ell}}\#\pi}(\rho)$.
Next consider the condition (\ref{eq:notationsystem.12}).
We have $\forall i(K(\nu_{i})\subset\mathcal{H}_{\max K(\nu_{i})}[\Tht])$ by (\ref{eq:notation.12rfl}), and
hence $\forall i(K(\nu_{i})\subset\mathcal{H}_{\max K(\nu_{i})}(\rho))$ by $\Theta\subset\rho$.
Therefore $\rho\in OT$ by Definition \ref{df:notationsystem}.\ref{df:notationsystem.12}.
Moreover $\rho\in C(\pi,\gam,\Tht)$, i.e., $\calh_{\gam}(\rho)\cap\pi\subset\rho\spand \Tht\cap\pi\subset\rho$.
Hence $\rho\in H(\vec{\nu},\pi,\gam,\Tht)$.



By Inversion lemma \ref{lem:inversionreg} we obtain for each $\del\equiv(\exi x\in L_{\pi}\del_{1}(x))\in\Del$ 
and each $d\in Tm(\rho)$ with $|d|=\max(\{0\}\cup {\sf k}(d))$,
$(\mathcal{H}_{\gam\# |d|},\Theta\cup{\sf k}(d))\vdash^{a_{\ell}}_{\pi}\Gamma,\lnot\del_{1}(d)$.

We have $\{\pi,\gam,|d|\}\subset\calh_{\gam\#|d|}(\pi)$ by $|d|<\rho<\pi$, and this yields
$|d|\in\calh_{\gam\#|d|}(\psi_{\pi}(\gam\#|d|))\cap\pi\subset\psi_{\pi}(\gam\#|d|)$.
Hence $|d|<\psi_{\pi}(\gam\#|d|)$, and
$\fal e(\Theta\cup{\sf k}(d)\subset\psi_{\pi}(\gam\#|d|\# e))$, i.e., (\ref{eq:asscollasKN}) holds for $\lam=\pi$
and $\gam\#|d|$.
Let $\bet_{d}=\psi_{\pi}(\widehat{a_{d}})$ for $\widehat{a_{d}}=\gam\#|d|\#\ome^{\pi+a_{\ell}+1}=\widehat{a_{\ell}}\#|d|$.
SIH yields 
$(\calh_{\widehat{a_{d}}+1},\Tht\cup{\sf k}(d))\vdash^{\bet_{d}}_{\bet_{d}}\Gam,\lnot\del_{1}(d)$,
which in turn 
Boundedness lemma \ref{lem:boundednessreg} yields
$(\calh_{\widehat{a_{\pi}}+1},\Tht\cup{\sf k}(d))\vdash^{\bet_{d}}_{\bet_{d}}\Gam,\lnot\del_{1}^{(\bet_{d},\pi)}(d)$
 for 
$\widehat{a_{\pi}}=\gam\#\pi\#\ome^{\pi+a_{\ell}+1}=\widehat{a_{\ell}}\#\pi$.
By persistency we obtain 
$(\calh_{\widehat{a_{\pi}}+1},\Tht\cup{\sf k}(d))\vdash^{\bet_{d}}_{\rho}\Gam,\lnot\del_{1}^{(\rho,\pi)}(d)$
for $\bet_{d}<\psi_{\pi}(\widehat{a_{\pi}})=\rho\in\calh_{\gam}[\Tht]$.
Since $d\in Tm(\rho)$ is arbitrary, $(\bigwedge)$ yields 
\beqn\label{eq:collapsKNcase2}
(\calh_{\widehat{a_{\pi}}+1},\Tht)\vdash^{\rho}_{\rho}\Gam,\lnot\del^{(\rho,\pi)}
\eeqn

Now pick the $\rho$-th branch from the right upper sequents
\[
 (\mathcal{H}_{\widehat{a_{\pi}}+1},\Theta\cup\{\rho\}\vdash^{a_{r}(\rho)}_{\pi}
 \Gamma, \Del^{(\rho,\pi)}
 \]
By $\rho\in\mathcal{H}_{\widehat{a_{\pi}}+1}[\Theta]$ and Lemma \ref{lem:weakening}.\ref{lem:weakening2} we obtain
\beqn\label{eq:collapsKNcase2a}
 (\mathcal{H}_{\widehat{a_{\pi}}+1},\Theta)\vdash^{a_{r}(\rho)}_{\pi}
 \Gamma, \Del^{(\rho,\pi)}
\eeqn
{\bf Case 2.1}.
First consider the case $\lam=\pi$. Then $\Del^{(\rho,\pi)}\subset\Sig_{0}(\lam)$.
Let $\bet_{\rho}=\psi_{\pi}(b_{\rho})$ with 
$b_{\rho}=\widehat{a_{\pi}}\#1\#\ome^{\pi+a_{r}(\rho)+1}=\gam\#\ome^{\pi+a_{\ell}+1}\#\ome^{\pi+a_{r}(\rho)+1}\#\pi\#1$.
Then $\bet_{\rho}>\rho$ and $\fal d[\Tht\cup\{\rho\}\subset\psi_{\pi}(\widehat{a_{\pi}}+1\#d)]$.
SIH yields for (\ref{eq:collapsKNcase2a}) 
\beqn\label{eq:collapsKNcase2b}
 (\mathcal{H}_{b_{\rho}+1},\Theta)\vdash^{\bet_{\rho}}_{\bet_{\rho}}
 \Gamma, \Del^{(\rho,\pi)}
 \eeqn
Several $(cut)$'s with (\ref{eq:collapsKNcase2b}), (\ref{eq:collapsKNcase2})
yield
$(\calh_{\hat{a}+1},\Tht)\vdash^{\bet_{\rho}+p}_{\bet_{\rho}} \Gam$
for $\bet_{\rho}\geq\rho$, $\widehat{a_{\pi}}<b_{\rho}<\hat{a}$ and some $p<\ome$, 
where $\bet_{\rho}<\bet=\psi_{\pi}(\hat{a})$ 
by $b_{\rho}<\hat{a}$.
(\ref{eq::collapsKNconcl}) follows.
\\

\noindent
{\bf Case 2.2}.
Next consider the case when $\lam<\pi$.
Then $\lam<\rho$ and $\Del^{(\rho,\pi)}\subset\Sig_{1}(\rho^{+})$ with $\rho^{+}=\Ome_{\rho+1}$.
SIH with (\ref{eq:collapsKNcase2a}) yields 
$ (\mathcal{H}_{b_{\rho}+1},\Theta\cup\{\rho\})\vdash^{\bet_{\rho^{+}}}_{\bet_{\rho^{+}}}
 \Gamma, \Del^{(\rho,\pi)}$
 for $\bet_{\rho^{+}}=\psi_{\rho^{+}}(b_{\rho})>\rho$,
and by Lemma \ref{lem:weakening}.\ref{lem:weakening2} we obtain
 \beqn\label{eq:collapsKNcase2d}
 (\mathcal{H}_{b_{\rho}+1},\Theta)\vdash^{\bet_{\rho^{+}}}_{\bet_{\rho^{+}}}
 \Gamma, \Del^{(\rho,\pi)}
 \eeqn
Several $(cut)$'s with (\ref{eq:collapsKNcase2}), (\ref{eq:collapsKNcase2d}) 
yield
$(\calh_{b_{0}+1},\Tht)\vdash^{\bet_{\rho^{+}}+p}_{\bet_{\rho^{+}}} \Gam$
for
$\bet_{\rho^{+}}>\rho$ and
$b_{0}=\gam\#(\ome^{\pi+a_{\ell}+1}\cdot 2)\#\ome^{\pi+a_{r}(\rho)+1}\#1\geq\max\{b_{\ell},b_{\rho}\}$.
Predicative cut-elimination lemma \ref{lem:predcereg} yields 
for $\bet_{1}=\vphi(\bet_{\rho^{+}})(\bet_{\rho^{+}}+p)<\rho^{+}$
\beqn\label{eq:mu00}
(\calh_{b_{0}+1},\Tht)\vdash^{\bet_{1}}_{\rho} \Gam
\eeqn
We obtain
 $\lam<\rho\in\calh_{b_{0}+1}[\Tht]$
by $\gam<\widehat{a_{\ell}}< b_{0}$.
 MIH with (\ref{eq:mu00}) 
yields
$(\calh_{c+1},\Tht)\vdash^{\psi_{\lam}c}_{\psi_{\lam}c} \Gam$ for $c=b_{0}\#1\#\ome^{\rho+\bet_{1}+1}$.
We obtain
$c=b_{0}\#\ome^{\rho+\bet_{1}+1}\#1=\gam\#(\ome^{\pi+a_{\ell}+1}\cdot 2)\#\ome^{\pi+a_{r}(\rho)+1}\#\ome^{\rho+\bet_{1}+1}\#2
<\gam\#\ome^{\pi+a+1}=\hat{a}$
since 
$a_{\ell}, a_{r}(\rho)<a$ and $\rho,\bet_{1}<\rho^{+}<\pi$.
Hence
 $\psi_{\lam}c<\psi_{\lam}(\hat{a})=\bet$, and (\ref{eq::collapsKNconcl}) follows.
 \\
 
 \noindent
{\bf Case 3}. 
Third consider the case when the last inference introduces a $\Sig_{1}(\lam)$-sentence 
$(\fal x\in c\,\tht(x))\in\Gam$ for $c\in Tm(\lam)$.
We have 
$(\mathcal{H}_{\gam},\Theta\cup{\sf k}(d))\vdash^{a(d)}_{\pi}\Gamma,\tht(d)$ for each
$d\in Tm(|c|)$.
Then we see from (\ref{eq:collapsKNclose}) that
$|d|<|c|\in\calh_{\gam}(\psi_{\rho}(\gam\#e))\cap\rho\subset\psi_{\rho}(\gam\#e)$ for 
 any $\rho\in[\lam,\pi]$ and any $e$.
 Hence $|d|\in\psi_{\rho}(\gam\#e)$.
(\ref{eq:asscollasKN}) is enjoyed for $\Tht\cup{\sf k}(d)$.
SIH yields 
$
(\mathcal{H}_{\hat{a}+1},\Theta\cup{\sf k}(d))\vdash^{\bet_{d}}_{\bet_{d}}\Gamma,
\tht(d)
$
for $\bet_{d}=\psi_{\lam}(\widehat{a(d)})$.
$(\bigwedge)$ yields (\ref{eq::collapsKNconcl}) for $\bet=\psi_{\lam}(\hat{a})>\bet_{d}$.
\\

\noindent
{\bf Case 4}.
Fourth consider the case when the last inference introduces a $\Sig_{1}(\lam)$-sentence 
$(\exi x\in L_{\lam}\,\tht(x))\in\Gam$.
We have $(\calh_{\gam},\Tht)\vdash^{a_{0}}_{\pi}\Gam,\tht(d)$ for a $d\in Tm(\lam)$.
SIH yields $(\calh_{\hat{a}+1},\Tht)\vdash^{\bet_{0}}_{\bet_{0}}\Gam,\tht(d)$
for $\bet=\psi_{\lam}(\hat{a})>\psi_{\lam}(\widehat{a_{0}})=\bet_{0}$.
Without loss of generality we can assume that ${\sf k}(d)\subset{\sf k}(\tht(d))$.
Then we see from (\ref{eq:collapsKNclose}) 
that 
$|d|\in\calh_{\gam}(\psi_{\lam}(\gam+1))\cap\lam\subset\psi_{\lam}(\gam+1)<\bet$.
Thus  is enjoyed in the following inference rule $(\bigvee)$.
We obtain $(\calh_{\hat{a}+1},\Tht)\vdash^{\bet}_{\bet}\Gam$
by a $(\bigvee)$, which enjoys (\ref{eq:bigveebnd}).
\\

\noindent
{\bf Case 5}.
Fifth consider the case when the last inference is a $({\rm rfl}(\tau, j, \vec{\mu},\vec{\nu}))$
for a $\tau\in\calh_{\gam}[\Tht]\cap\pi$.
We have an $a_{\ell}<a$ and a finite set $\Del$ of $\Sig_{j}(\tau)$-sentences such that 
$(\mathcal{H}_{\gam},\Tht)\vdash^{a_{\ell}}_{\pi}\Gamma,\lnot\del$
for each $\del\in\Del$.
On the other hand we have a sequence $\vec{\nu}$ and
an ordinal term $a_{r}(\rho)<a$ for each $\rho\in H(\vec{\nu},\tau,\gam,\Tht)$
such that
$(\mathcal{H}_{\gam},\Tht\cup\{\rho\})\vdash^{a_{r}(\rho)}_{\pi}\Gamma, \Del^{(\rho,\tau)}$.
By (\ref{eq:collapsKNclose}),
for any $\rho\in H(\vec{\nu},\tau,\gam,\Tht)$ we obtain
\beqn\label{eq:case5}
\fal e\fal \kap[\max\{\tau+1,\lam\}\leq\kap\leq\pi \Rarw \rho<\tau\in\calh_{\gam}(\psi_{\kap}(\gam\#e))\cap\kap\subset\psi_{\kap}(\gam\#e)]
\eeqn
{\bf Case 5.1}.
First consider the case when $\tau<\lam$.
Then $\rho<\psi_{\kap}(\gam\#e)$ for any $\kap\in[\lam,\pi]$ and $e$.
From SIH with (\ref{eq:case5})
we obtain 
$(\mathcal{H}_{\hat{a}+1},\Tht)\vdash^{\bet_{\ell}}_{\bet_{\ell}}\Gamma,\lnot\del$
for each $\del\in\Del$ with $\bet_{\ell}=\psi_{\lam}(\widehat{a_{\ell}})$, and
$(\mathcal{H}_{\hat{a}+1},\Tht\cup\{\rho\} )\vdash^{\bet_{r}(\rho)}_{\bet_{r}(\rho)}\Gamma, \Del^{(\rho,\tau)}$
for each $\rho\in H(\vec{\nu},\tau,\gam,\Tht)$ with $\bet_{r}(\rho)=\psi_{\lam}(\widehat{a_{r}(\rho)})$.
We see $\max\{\bet_{\ell},\bet_{r}(\rho),\tau\}<\bet=\psi_{\lam}(\hat{a})$, and
an inference rule $({\rm rfl}(\tau,j,\vec{\mu},\vec{\nu}))$
yields $(\mathcal{H}_{\hat{a}+1},\Tht)\vdash^{\bet}_{\bet}\Gam$.
\\

\noindent
{\bf Case 5.2}.
Second consider the case when $\lam\leq\tau$.
Then $\Del\cup\Del^{(\rho,\tau)}\subset\Sig_{1}(\tau^{+})$, and
$\rho<\psi_{\kap}(\gam\#e)$ for $\tau<\kap\leq\pi$ and $e$ by (\ref{eq:case5}).
SIH yields 
$(\mathcal{H}_{\widehat{a_{\ell}}+1},\Tht)\vdash^{\bet_{2}}_{\bet_{2}}\Gamma,\lnot\del$
for each $\del\in\Del$, where
$\bet_{2}=\psi_{\tau^{+}}\left(\widehat{a_{\ell}}\right)$.
On the other side SIH yields
$(\mathcal{H}_{\widehat{a_{r}(\rho)}+1},\Tht\cup\{\rho\} ) \vdash^{\bet_{\rho}}_{\bet_{\rho}}
 \Gamma, \Del^{(\rho,\tau)}$
 for each $\rho\in H(\vec{\nu},\tau,\gam,\Tht)$, where $\bet_{\rho}=\psi_{\tau^{+}}\left(\widehat{a_{r}(\rho)}\right)$.
Predicative cut-elimination lemma \ref{lem:predcereg} yields 
$(\mathcal{H}_{\widehat{a_{\ell}}+1},\Tht)\vdash^{\del_{2}}_{\tau}\Gamma,\lnot\del$
and
$(\mathcal{H}_{\widehat{a_{r}(\rho)}+1},\Tht\cup\{\rho\})\vdash^{\del_{\rho}}_{\tau}
 \Gamma, \Del^{(\rho,\tau)}$
for
$\del_{2}=\vphi(\bet_{2})(\bet_{2})$ and $\del_{\rho}=\vphi(\bet_{\rho})(\bet_{\rho})$.
From these with the inference rule $({\rm rfl}(\tau,j,\vec{\mu},\vec{\nu}))$ we obtain
\beqn\label{eq:tauxi}
(\calh_{\widehat{a_{0}}+1},\Tht)\vdash_{\tau}^{\del_{0}+1}\Gam
\eeqn
where $\sup\{\del_{2},\del_{\rho}: \rho\in H(\vec{\nu},\tau,\widehat{a_{0}}+1,\Tht)\}\leq\del_{0}:=\vphi(\bet_{0})(\bet_{0})\in\calh_{\widehat{a_{0}}+1}[\Tht]$
with 
$\sup\{\bet_{2},\bet_{\rho}: \rho\in H(\vec{\nu},\tau,\gam,\Tht)\}\leq\bet_{0}:=\psi_{\tau^{+}}\left(\widehat{a_{0}}\right)$,
and
$\sup\{a_{\ell},a_{r}(\rho): \rho\in H(\vec{\nu},\tau,\gam,\Tht)\}\leq a_{0}\in\calh_{\gam}[\Tht]\cap a$, cf.\,(\ref{eq:Mhordgap}).


MIH with (\ref{eq:tauxi}) 
yields $(\calh_{\hat{a}+1},\Tht)\vdash_{\del}^{\del}\Gam$ 
for $\del=\psi_{\lam}((\widehat{a_{0}}+1)\#\ome^{\tau+\del_{0}+2})$
and $(\widehat{a_{0}}+1)\#\ome^{\tau+\del_{0}+2}<\hat{a}$.
We have $\del=\psi_{\lam}(\widehat{a_{0}}\#1\#\ome^{\tau+\del_{0}+2})<\psi_{\lam}(\hat{a})=\bet$
by $\widehat{a_{0}}<\hat{a}$ and $\tau,\del_{0}<\tau^{+}<\pi$ and $\tau\in\calh_{\gam}[\Tht]$.
(\ref{eq::collapsKNconcl}) follows.
\\

\noindent
{\bf Case 6}.
Sixth consider the case when the last inference is a $(cut)$.
For an $a_{0}<a$ and a $C$ with $\mbox{rk}(C)<\pi$, we have
$(\calh_{\gam},\Tht)\vdash^{a_{0}}_{\pi}\Gam,\lnot C$ and
$(\calh_{\gam},\Tht)\vdash^{a_{0}}_{\pi}C,\Gam$.
\\

\noindent
{\bf Case 6.1}. First consider the case when $\rk(C)<\lam$.
Then $C\in\Sig_{0}(\lam)$.
SIH yields the lemma.
\\

\noindent
{\bf Case 6.2}. Second consider the case when $\lam\leq\rk(C)<\pi$.
Let $\rho^{+}=(\rk(C))^{+}=\min\{\kap\in Reg: \rk(C)<\kap\}$. 
Then $C\in\Sig_{0}(\rho^{+})$ and 
$\lam\leq\rho\in\calh_{\gam}[\Tht]\cap\pi$.
SIH yields 
$(\calh_{\widehat{a_{0}}+1},\Tht)\vdash^{\bet_{0}}_{\bet_{0}}\Gam,\lnot C$
and
$(\calh_{\widehat{a_{0}}+1},\Tht)\vdash^{\bet_{0}}_{\bet_{0}}C,\Gam$
for $\bet_{0}=\psi_{\rho^{+}}\left(\widehat{a_{0}}\right)\in\calh_{\widehat{a_{0}}+1}[\Tht]$.
By a $(cut)$ we obtain
$(\calh_{\widehat{a_{0}}+1},\Tht)\vdash^{\bet_{1}}_{\bet_{1}}\Gam$
 for $\bet_{1}=\max\{\bet_{0},\rk(C)\}+1$ with $\rho<\bet_{1}<\rho^{+}$.
Predicative cut-elimination lemma \ref{lem:predcereg} yields 
$(\calh_{\widehat{a_{0}}+1},\Tht)\vdash^{\del_{1}}_{\rho}\Gam$
for $\del_{1}=\vphi(\bet_{1})(\bet_{1})$,
where $\widehat{a_{0}}\in\calh_{\widehat{a_{0}}+1}[\Tht]$, and 
$\fal e\fal \tau\in[\lam,\rho][\Tht\subset\psi_{\tau}(\widehat{a_{0}}\#e)]$ hold.
Hence MIH with $\rho\in\calh_{\widehat{a_{0}}+1}[\Tht]$ 
yields $(\calh_{b+1},\Tht)\vdash^{\psi_{\lam}(b)}_{\psi_{\lam}(b)}\Gam$
for $b=\widehat{a_{0}}\#1\#\ome^{\rho+\del_{1}+1}$.
We see $b<\hat{a}$ and $\psi_{\lam}(b)<\psi_{\lam}(\hat{a})=\bet$, and (\ref{eq::collapsKNconcl}) follows.
\\

\noindent
{\bf Case 7}.
Seventh consider the case when the last inference is an $(\Ome\in M_{2})$.
We have $(\calh_{\gam},\Tht)\vdash^{a_{\ell}}_{\pi}\Gam,C$ for an $a_{\ell}<a$, and
$(\calh_{\gam},\Tht\cup\{\alp\})\vdash^{a_{r}(\alp)}_{\pi}\lnot C^{(\alp,\Ome)}, \Gam$
for an $a_{r}(\alp)<a$ for each $\alp<\Ome$,
where $C\in\Pi_{2}(\Ome)$.

The case $\lam>\Ome$ is seen as in {\bf Case 5.1}.
The case $\lam=\Ome$ is seen as in {\bf Case 5.2}.
\eprf
\\

Let us conclude Theorem \ref{thm:2}.
Let $\Ome=\Ome_{1}$.
\\

\noindent
{\bf Proof} of Theorem \ref{thm:2}.
Let ${\sf KP}\Pi_{N}\vdash\tht$.
By Embedding lemma \ref{th:embedregthm} pick an $m$ so that 
 $(\calh_{0},\emptyset)\vdash_{\mK+m}^{\mK\cdot 2+m}\tht$.
Predicative cut-elimination lemma \ref{lem:predcereg} yields 
$(\calh_{0},\emptyset)\vdash_{\mK}^{\ome_{m+1}(\mK+1)}\tht$
for $\ome_{m}(\mK\cdot 2+m)<\ome_{m+1}(\mK+1)$.
Lemma \ref{lem:collapsKN} yields 
$(\calh_{a+1},\emptyset)\vdash^{\bet}_{\bet}\tht$
for $a=\ome^{\mK+\ome_{m+1}(\mK+1)+1}$ and $\bet=\psi_{\Ome}(a)$.
Predicative cut-elimination lemma \ref{lem:predcereg} yields
$(\calh_{a+1},\emptyset)\vdash^{\vphi(\bet)(\bet)}_{0}\tht$.
We obtain $\vphi(\bet)(\bet)<\alp:=\psi_{\Ome}(\ome_{n}(\mK+1))$ for $n=m+3$, and hence
$(\calh_{\ome_{n}(\mK+1)},\emptyset)\vdash^{\alp}_{0}\tht$.
Boundedness lemma \ref{lem:boundednessreg} yields
$(\calh_{\ome_{n}(\mK+1)},\emptyset)\vdash^{\alp}_{0}\tht^{(\alp,\Ome)}$.
Since each inference rule other than reflection rules $({\rm rfl}(\pi,k,\vec{\xi},\vec{\nu}))$ and $(\Omega\in M_{2})$ is sound,
we see by induction up to $\alp=\psi_{\Ome}(\ome_{n}(\mK+1))$ that
$L_{\alp}\models\tht$.

This completes a proof of Theorem \ref{thm:2}.


\begin{thebibliography}{99}









\bibitem{LMPS}
T. Arai, 
Iterating the recursively Mahlo operations,
in \textit{Proceedings of the thirteenth International Congress of Logic Methodology, Philosophy of Science}, 
eds. C. Glymour, W. Wei and D. Westerstahl
(College Publications, King's College, London, 2009), pp. 21--35.

\bibitem{WFnonmon2}
T. Arai, 
Wellfoundedness proofs by means of non-monotonic inductive definitions II: first order operators,
\textit{Ann. Pure Appl. Logic} {\bf 162} (2010), pp. 107-143.




\bibitem{consv}T. Arai, 
Conservations of first-order reflections,
 \textit{Jour. Symb. Logic} {\bf 79} (2014), pp. 814-825.


\bibitem{PTPiN}
T. Arai,  
Proof theory for theories of ordinals III:$\Pi_{N}$-reflection, 
\textit{Gentzen's centenary: the quest of consistency}, eds. R. Kahle and M. Rathjen,
pp. 357-424,
Springer, 2015.



\bibitem{KPpiNwfprf}T. Arai,
Wellfoundedness proof for first-order reflection,
draft.
posted to the arxiv.





\bibitem{Buchholz} W. Buchholz, 
A simplified version of local predicativity, 
in \textit{Proof Theory},
eds. P. H. G. Aczel, H. Simmons and S. S. Wainer
(Cambridge UP,1992), pp. 115--147.



\bibitem{Pohlers} W. Pohlers and J.-C. Stegert,
Provable recursive functions of reflection,
in \textit{Logic, Construction, Computation}, eds. U. Berger, H. Diener, P. Schuster and M. Seisenberger,
Ontos Mathematical Logic vol. 3 (De Gruyter, 2012), pp. 381-474.


\bibitem{Rathjen94} M. Rathjen, 
Proof theory of reflection, 
\textit{Ann. Pure Appl. Logic} {\bf 68} (1994) 181--224.

\bibitem{RathjenAFML1}M. Rathjen, 
An ordinal analysis of stability, 
\textit{Arch. Math. Logic} {\bf 44} (2005) 1-62.



\end{thebibliography}
\end{document}